\documentclass{amsart}
\usepackage{graphicx}
\usepackage{latexsym}
\usepackage{epsf}
\usepackage{epsfig}
\usepackage{color}
\usepackage{a4}
\usepackage{amsmath} 
\usepackage{amssymb}
\usepackage{amsxtra}
\usepackage{fancyhdr}
\usepackage{fullpage}
\usepackage{enumitem}
\usepackage{caption}
\usepackage{subcaption}
\usepackage{tikz}
\usepackage{tikz-3dplot}
\usepgflibrary{arrows}
\usepackage[utf8]{inputenc}
\usepackage[english]{babel}
\usepackage{lineno,hyperref}
\modulolinenumbers[5]

\newtheorem{lemma}{Lemma}
\newtheorem{proposition}{Proposition}
\newtheorem{corollary}{Corollary}
\newtheorem{theorem}{Theorem}
\newtheorem{remark}{Remark}
\newtheorem{definition}{Definition}

\newcommand{\OO}{\mathcal {O}}
\newcommand{\Q}{\mathcal{Q}}
\renewcommand{\P}{\mathcal{P}}

\newcommand{\R}{\mathcal{R}}
\newcommand{\T}{\mathcal{T}}
\newcommand{\bo}{\rm bo}

\newcommand{\mfT}{\mathfrak T}
\newcommand{\eps}{\epsilon}
\newcommand{\dsp}{\displaystyle}
\newcommand{\RR}{\mathbb{R}}
\newcommand{\NN}{\mathbb{N}}
\newcommand{\mfQ}{\mathfrak{Q}}
\newcommand{\esssup}{\displaystyle{{\rm ess\,sup}}}
\newcommand{\zetabar}{\underline{\zeta}}

\newcommand{\vbar}{\underline{v}}

\newcommand{\nn}{\nonumber}

\newcommand{\defeq}{\overset{\text{\tiny def}}{=}}

\makeindex

\begin{document}

\title{Medium amplitude model for internal waves over large topography variation}

%
\author{Ralph Lteif}
\address{Lebanese American University (LAU), Graduate Studies and Research (GSR), School of Arts and Sciences, Computer Science and Mathematics Department, Beirut, Lebanon}
\email{ralph.lteif@lau.edu.lb}
\author{Bashar Khorbatly}
\address{Lebanese American University (LAU), Graduate Studies and Research (GSR), School of Arts and Sciences, Computer Science and Mathematics Department, Beirut, Lebanon}
\email{bashar-elkhorbatly@hotmail.com}

\date{ \today}

\begin{abstract}
The purpose of this paper is to present the derivation and mathematical analysis of a new asymptotic model that describes the evolution of medium amplitude internal waves propagating between a flat rigid-lid and a highly variable topography. The smallness assumptions on the topographic variation parameter used in [Communications on Pure \& Applied Analysis, 2015, 14 (6): 2203-2230] and in [Asymptotic Analysis, vol. 106, no. 2, pp. 61-98, 2018] are now relaxed and the results of the aforementioned papers are improved and generalized to the complex case of large topography variation. Limiting the flow to one-layer, we also emphasize our  model's well-posedness in comparison to the original asymptotic model.
\end{abstract}

\maketitle

\tableofcontents

\section{Introduction}

\subsection{Motivation}
Ocean water is not uniform when it comes to mass density. In fact, the temperature and salinity of water in the ocean vary according to depth which create a stratification effect dividing the water into layers with different densities. The disturbance of these layers by tidal flows over variable topography generates internal waves. The absorbed solar radiation make the upper surface water warmer with a lower density lying above a colder denser water. Internal waves play an important role in underwater biological life and navigation, thus understanding their behavior is very essential. Herein, we consider the uni-dimensional flow of internal waves. Simplifying assumptions on the nature of the fluids are commonly used in oceanography. Namely, the fluids are supposed to be homogeneous, immiscible, inviscid and affected only by gravitational force. Moreover, irrotational and incompressible fluids are considered.

The mathematical aspects of the internal waves have been the subject of many studies in the literature. The internal wave flow is described by two evolution equations. These equations are commonly called the ``\emph {full Euler system}". We omit here the detailed derivation of the Zakharov formulation of this system for the sake of readability, instead, the equations governing the two-layer flow are briefly recalled in Section~\ref{FEsec}. For more details, the interested reader can see for instance~\cite{Anh09,BonaLannesSaut08,Duchene13}. Solving the ``\emph {full Euler system}" mathematically is very difficult due to the free boundary problem \emph{i.e} the domain is moving with time. Indeed, the interface deformation boundary function is one of the unknowns. To overcome this problem, many researchers searched for approximate solutions of the exact system. To this end, the derivation of simpler asymptotic models had attract a lot of attention. These approximate models are derived in much simpler settings where dimensionless variables and unknowns are introduced, allowing a fair description of the exact behavior of the full system in particular physical regimes. In this study, we construct a new asymptotic model for the two-layer flow over strong variations in bottom topography using an additional smallness assumption on the interface deformation. More precisely, we consider medium amplitude deformations at the interface level. This smallness assumption corresponds to the commonly known Camassa-Holm (CH) regime.

The two-layer flow has been widely studied in the literature, laying the groundwork for an important theoretical framework. Many approximate models describing the evolution of a two-layer flow over flat topography and under a rigid lid have previously been developed and studied, see for instance~\cite{Miyata85,Matsuno93,ChoiCamassa96,ChoiCamassa99,BonaLannesSaut08,GuyenneLannesSaut10,DucheneIsrawiTalhouk15} and references therein. Two-layer flow over variable topography has been also derived and studied in some significant works~\cite{NachbinZarate09,Anh09,Duchene10,ChoiBarros13,
DucheneIsrawiTalhouk14}. These models were proved to be consistent with the exact system, however they are not supported with a full justification (\emph{i.e} well posedness, consistency with respect to the full Euler system and convergence to the latter). More recently, in an effort to deal with medium amplitude internal waves propagating over medium amplitude topography variations, a full justification result for a newly derived asymptotic model of Green-Nadghi type in the Camassa-Holm regime (GNCH) has been obtained in~\cite{IsrawiLteifTalhouk15}. In~\cite{LteifIsrawi18}, the full justification result of the GNCH model has been improved to a more complex case of large amplitude topography with slow variation.
In this paper, we present the derivation and mathematical analysis of a new GNCH model allowing for strong topography variations. This new model presents an important improvement upon the existing ones introduced in~\cite{IsrawiLteifTalhouk15,LteifIsrawi18} since both smallness and slow variation assumptions on bottom topography are now relaxed which is quite more rational in ocean floor. Nevertheless, some serious difficulties arise. Indeed, relaxing any smallness assumption on the amplitude parameter of topographic variation prompt some terms that are actually accompanied with derivative terms on the interface deformation that cannot be controlled by the intended energy norm associated with our model and thus do not always allow its full justification. To overcome this difficulty, we establish a specially designed model allowing to deal with these terms and possessing a quasilinear hyperbolic structure.
Classical theory of hyperbolic systems is then applied, in particular energy estimates, hence allowing its well-posedness. Moreover, the full justification of the newly derived model follows in the same way as the models studied in   \cite{DucheneIsrawiTalhouk14,IsrawiLteifTalhouk15,LteifIsrawi18}. The obtained model is valid under some key restrictions on the bottom deformation, see Remark~\ref{oderem}. Moreover, restricting the flow to one-layer, we stress our model's well-posedness in comparison to the original asymptotic model studied in~\cite{Haidar2018}.

\subsection{Organization of the paper}
In Section \ref{FEsec}, we briefly introduce the \emph{full Euler system}. In Section \ref{Sec3}, we precisely derive the new asymptotic model starting from the original Green-Naghdi model. Some preliminary results on the properties of the specially designed symmetric differential operator are stated in Section \ref{Sec4}. Section~\ref{The-linearized-system} contains an essential procedure preceding the proof of our result. Section~\ref{Sec6} is devoted to the linear analysis of the asymptotic model and its well-posedness result. In Section \ref{Sec7}, we limit the work to the one-layer case and emphasize the well-posedness of our model in comparison to the original asymptotic model.

\subsection{Notations}
We refer to $C_0$ as a nonnegative constant whose exact expression isn't important. Consider the notation $ a \lesssim b$ as $a\leq C_0\ b$. Also consider $A=\OO(B)$ as $A\leq C_0\ B$.\\
$C(\lambda_1, \lambda_2,..)$ denotes a nonnegative constant that depends on the parameters $\lambda_1$, $\lambda_2$.. and whose dependence on the $\lambda_j$ is always assumed to be nondecreasing.\\
$L^p=L^p(\RR)$ is the space of all Lebesgue-measurable functions
 $f$ with $1\leq p< \infty$ endowed with the norm $ \vert \psi \vert_{L^p}=\big(\int_{\RR}\vert \psi(x)\vert^p dx\big)^{1/p}<\infty$.  
If $p=\infty$, $L^\infty=L^\infty(\RR)$ consists of all essentially bounded, Lebesgue-measurable functions
 $f$ with the norm
$
 \vert \psi\vert_{L^\infty}= \esssup \vert \psi(x)\vert<\infty.
$\\
 Let $k \in \NN$, we denote by $W^{k,\infty}=W^{k,\infty}(\RR)=\{f \in L^{\infty}, |f|_{W^{k,\infty}}< \infty\}$, where $ |f|_{W^{k,\infty}}=\displaystyle{\sum_{\alpha \in \NN , \alpha\leq k} |\partial_x^{\alpha} f|_{L^{\infty}}}$.
Denote by $(1-\partial_x^2)^{1/2}$ the pseudo-differential operator $\Lambda$. For any positive number $s$, $H^s=H^s(\RR)$ stands for the Sobolev space of all tempered
 distributions $f$ with norm $\vert f\vert_{H^s}=\vert \Lambda^s f\vert_{L^2} < + \infty$.
 We refer to $L^\infty([0,T);H^s(\RR))$ the space of functions such that $u(t,\cdot)$ is controlled in $H^s$, uniformly for $t\in[0,T)$:
\begin{equation*} \big\Vert u\big\Vert_{L^\infty([0,T);H^s(\RR))} \ = \ \esssup_{t\in[0,T)} \vert u(t,\cdot)\vert_{H^s} \ < \ \infty.\end{equation*} 
Let $L$ be any closed operator defined and on a Banach space $Z$. The commutator $[L,\psi]\Psi=L(\psi\Psi )-\psi L(\Psi)$  is defined for any $\psi$, $\Psi$, and $\psi\Psi$ belonging to $L$'s domain. The same notation is used for $\psi$, an operator that maps $L$'s domain into itself.
\section{Full Euler system}\label{FEsec}
\begin{figure}[h]
\centering 
\includegraphics[scale=0.8]{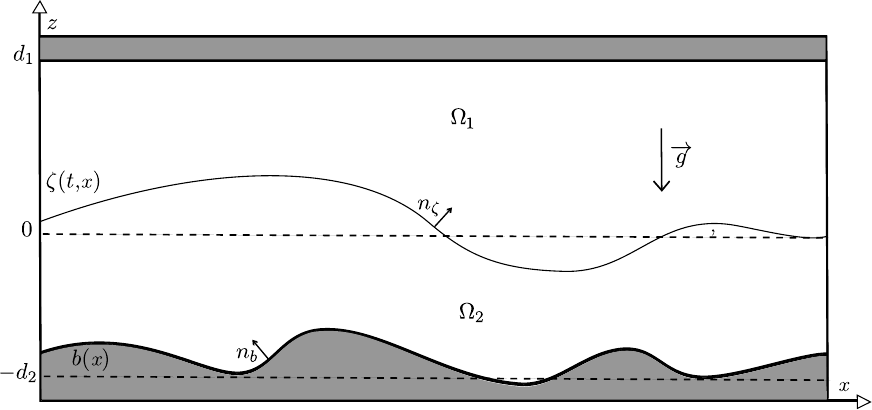}
\caption{Domain of study}
\label{dom1}
\end{figure}

In this section, we briefly recall the governing equations of the two-layer flow. We restrict our study to the one-dimensional setting, we assume that the internal wave represented by the function $\zeta(t,x)$ propagates between a flat rigid lid located at $d_1$ and a variable topography whose deformation with respect to its rest place $(x,-d_2)$ is represented by the function $b(x)$. The upper and lower layer domains denoted by $ \Omega_1$ and $\Omega_2$ respectively are supposed to stay connected, in other words the heights of the upper and lower fluids ($h_1$, $h_2$ respectively) must remain strictly positive, that is $h_1,h_2\geq h_0>0$,  for some constant height $h_0$, see Figure~\ref{dom1}.

Now let us specify the assumptions on the nature and domain of the fluids. This type of reasonable hypothesis is commonly used in oceanography when determining the governing equations of two-layer flow.

First, we consider both fluids are homogeneous, therefore the mass density of the top and low fluids denoted by $\rho_1$ and $\rho_2$ respectively are constant.  Each layer of fluid 
is incompressible so the corresponding velocity field has a zero divergence. Assuming irrotational flows, there exists velocity potentials denoted by $\phi_i$ $(i=1,2)$ that satisfy the Laplace equation. Assuming ideal fluids with no viscosity, one obtains two Bernoulli equations. The surface, interface and bottom are all assumed to be bounding surfaces, that is to say no particle of fluid can cross the surface, interface or bottom. One may close the set of equations by assuming that the pressure is continuous at the interface. At this stage, one obtains:
\everymath{\displaystyle}
  \begin{equation}  \label{eqn:EulerComplet}
\left\{\begin{array}{ll}
         \Delta_{x,z} \phi_i \ = \ 0 & \mbox{ in }\Omega_i, \ i=1,2,\\
       \rho_i  \partial_t \phi_i+\frac{\rho_i}{2} |\nabla_{x,z} \phi_i|^2=-P_{i}-\rho_i gz & \mbox{ in }\Omega_i, \ i=1,2, \\
        \partial_n \phi_1=n_t .\nabla_{x,z} \phi_1=\partial_z \phi_1=0 & \mbox{ on } \{(x,z),z=d_1\}, \\
         \partial_t \zeta  \ = \ \sqrt{1+|\partial_x\zeta|^2}\partial_{n}\phi_1 \ = \ \sqrt{1+|\partial_x\zeta|^2}\partial_{n}\phi_2  & \mbox{ on } \{(x,z),z=\zeta(t,x)\},\\
         \partial_n\phi_2= n_b . \nabla_{x,z} \phi_2 \ = \ 0 &  \mbox{ on } \{(x,z),z=-d_2+b(x)\}, \\
       \lfloor P(t,x) \rfloor  = \sigma\partial_x \Big(\frac1{\sqrt{1+|\partial_x\zeta|^2}}\partial_x\zeta\Big) & \mbox{ on } \{(x,z),z=\zeta(t,x)\},
         \end{array}
\right.
\end{equation}
where $ \lfloor P(t,x) \rfloor \defeq \lim\limits_{\chi\to 0} \Big(  P(t,x,\zeta(t,x)+\chi)  -  P(t,x,\zeta(t,x)-\chi)  \Big)$ and $\partial_n=n.\nabla_{x,z}$ is the upward normal derivative in the direction of the vector $n$ under consideration. We denote by  $n_t=(0,1)^{T}$, $n_{\zeta}=\dfrac{1}{\sqrt{1+|\partial_x \zeta|^2}}(-\partial_x \zeta,1)^{T}$ and $n_{b}=\dfrac{1}{\sqrt{1+|\partial_x b|^2}}(-\partial_x b,1)^{T}$ the unit outward normal vectors at the top rigid surface, interface and bottom respectively. Here $\sigma$ denotes the surface (or interfacial) tension coefficient.

Studying both theoretical and numerical aspects of system~\eqref{eqn:EulerComplet} remain difficult as the domain is one of the unknowns. 
At this point, extra assumptions are made on some parameters with no dimensions to derive reduced asymptotic models and thus seek approximate solutions of the exact system. The parameters with no dimensions are the following:
\begin{align}\label{dimparam}
\mu \equiv \dfrac{d_1^2}{\lambda^2}, \quad \epsilon \equiv \dfrac{a}{d_1},\quad \beta \equiv \dfrac{a_b}{d_1},  \quad \gamma = \dfrac{\rho_1}{\rho_2}, \quad \delta \equiv  \dfrac{d_1}{d_2},  \quad {\rm bo^{-1}} = \frac{\sigma}{ g(\rho_2-\rho_1) d_1^2},
\end{align} 
where we denote by $a$ (resp. $a_b$) the maximal elevation of the internal wave (resp. bottom topography) and $\lambda$ the horizontal length scale of the wave at the interface. Written in its dimensionless form, system~\eqref{eqn:EulerComplet} can be reduced to two equations with two unknowns $(\zeta, \psi\equiv\phi_1(t,x,\zeta(t,x))$, see~\cite{Zakharov68,CraigSulem93}. For shortness sake, we do not recall here the Zakharov formulation of the \emph{full Euler system}. For more details, the interested reader could see for instance~\cite{Anh09,BonaLannesSaut08,Duchene13,DucheneIsrawiTalhouk14}.

The restrictions made on the dimensionless parameters defined in~\eqref{dimparam} describe the regimes under consideration. To complete this section, we state below the regime considered in this paper. 
\begin{definition}[Regime]\label{regi}
The model~\eqref{eq:Serre2} is valid under the following extra restrictions defining the Camassa-Holm regime,
\begin{align}\label{eqn:defRegimeCHmr}
\P_{\rm CH}  \equiv  \Big\{\ 0< \mu  \leq \mu_{\max}, \ 0 \leq  \epsilon   \leq \min(M \sqrt{\mu},1) ,\  0 \leq \beta\leq \beta_{\max}, \  \delta_{\min} \leq \delta \leq \delta_{\max}, \Big.\nn \\ \Big. \ 0 \leq \gamma < 1, \  \text{and} \ \bo_{\min} \leq \bo \leq \infty \Big\},
\end{align}
 with $0\leq \mu_{\max}, M, \beta_{\max}, \delta^{-1}_{\min}, \delta_{\max} , \bo^{-1}_{\min}, \nu_{0}^{-1} < \infty$ and $\bo_{\min} \gg 6$.\\
We proceed by denoting without difficulty
\[ M_{\rm SW} \ \equiv \ \max\big\{\mu_{\max},\delta_{\min}^{-1},\delta_{\max},\bo_{\min}^{-1},\beta_{\max}\big\}, \quad M_{\rm CH} \ \equiv \ \max\big\{M_{\rm SW},M,\nu_{0}^{-1}\big\}.\]

\end{definition}
\section{Asymptotic models for internal waves}\label{Sec3}
In this section, we develop a new asymptotic model for the propagation of medium amplitude internal waves over large topographic variations. We recall the original two-fluid Green-Naghdi system in subsection \ref{sec3.1}. In \ref{section3}, we derive the medium amplitude large topographic variation model under \ref{regi}. The derivation of the modified system to be studied is covered in section \ref{sec3.3}.
\subsection{The original Green-Naghdi two-fluid system}\label{sec3.1}
The well-known Green-Naghdi model is obtained by plugging the asymptotic expansions of the Dirichlet-Neumann operators given in~\cite{Duchene10,DucheneIsrawiTalhouk14} into the \emph{full Euler system}  and performing simple computations while disregarding all terms of order $\mu^2$. At this point, the well-known ``shear mean velocity" variable $v$, which connects the upper and lower layer depth averaged vertically integrated horizontal velocities $u_1$ and $u_2$, is introduced, see~\cite{DucheneIsrawiTalhouk14}: \begin{equation}\label{defv}v\equiv  \dfrac{u_2}{h_2} - \gamma  \dfrac{u_1}{h_1},\end{equation} 
$$ u_1   = \int_{\epsilon\zeta(t,x)}^{1} \frac{\partial}{\partial_x} \phi_1(t,x,z) \ dz\qquad\quad\text{ and }\qquad\quad u_2   =  \int_{-1/\delta+\beta b(x)}^{\epsilon\zeta(t,x)} \frac{\partial}{\partial_x} \phi_2(t,x,z) \ dz.$$
We do not define the Dirichlet-Neumann operators or their asymptotic expansions in this paper for the sake of simplicity. Instead, we briefly recall the Green-Naghdi system below and refer to~\cite{Duchene10,DucheneIsrawiTalhouk14,DucheneIsrawiTalhouk15} for more information on the derivation:
\begin{equation}\label{eq1}
\left\{ \begin{array}{l}
\displaystyle\partial_{ t}{\zeta} \ + \ \partial_x \Big(H(\epsilon\zeta,\beta b) v\Big)\ =\ 0,  \\ \\
\displaystyle\partial_{ t}\Big(  v \ + \ \mu\overline{\Q}[\epsilon\zeta,\beta b] v \Big) \ + \ (\gamma+\delta)\partial_x{\zeta} \ + \ \frac{\epsilon}{2} \partial_x\Big(H'(\epsilon\zeta,\beta b)  v^2\Big) \ = \displaystyle\mu\epsilon\partial_x\big(\overline{\R}[\epsilon
\zeta,\beta b] v  \big) +\mu\frac{\gamma+\delta}{\bo}\partial_{x}^3\zeta,
\end{array}
\right.
\end{equation}
where we denote by $H(\epsilon\zeta,\beta b)  = \dsp\frac{h_1h_2}{h_1+\gamma h_2}$, 
$H'(\epsilon\zeta,\beta b)  = \dsp\frac{h_1^2-\gamma h_2^2}{(h_1+\gamma h_2)^2}$, $h_1=1-\epsilon\zeta$ and $h_2= 1/\delta +\epsilon\zeta-\beta b$, as well as
\begin{align}\label{qu1}
\overline{\Q}[\epsilon\zeta,\beta b]v &=-\gamma\T[h_1,0]\Big(\frac{-h_2v}{h_1+\gamma h_2}\Big)+\T[h_2,\beta b]\Big(\frac{h_1v}{h_1+\gamma h_2}\Big) \nonumber\\
&=-\gamma\Big[\dfrac{1}{3h_1}\partial_x\Big(h_1^3\partial_x\big(\dfrac{h_2v}{h_1+\gamma h_2}\big)\Big)\Big]-\dfrac{1}{3h_2}\partial_x\Big(h_2^3\partial_x\big(\dfrac{h_1v}{h_1+\gamma h_2}\big)\Big)\nonumber\\ 
 &\quad +  \dfrac{1}{2h_2}\beta\Big[\partial_x\Big(h_2^2(\partial_x b) \dfrac{h_1 v}{h_1+\gamma h_2}\Big)-h_2^2(\partial_x b)\partial_x \big(\dfrac{h_1v}{h_1+\gamma h_2}\big)\Big]
+ \beta^2(\partial_x b)^2\big(\dfrac{h_1v}{h_1+\gamma h_2}\big)\;, 
\end{align}
\begin{align}\label{rbar}
\overline{\R}[\epsilon\zeta,\beta b] v & =-\dfrac{\gamma}{2}\Big(h_1\partial_x(\frac{-h_2v}{h_1+\gamma h_2})\Big)^2+ \dfrac{1}{2}\Big(-h_2\partial_x(\frac{h_1v}{h_1+\gamma h_2})+\beta(\partial_x{b})(\frac{h_1v}{h_1+\gamma h_2})\Big)^2\nonumber\\
& \quad+\gamma(\frac{-h_2v}{h_1+\gamma h_2})\T[h_1,0]\Big(\frac{-h_2v}{h_1+\gamma h_2}\Big)- (\frac{h_1v}{h_1+\gamma h_2})\T[h_2,\beta b]\Big(\frac{h_1v}{h_1+\gamma h_2}\Big)\; ,\qquad\quad\;
\end{align}
with $\T[h,b]\Psi\equiv \dfrac{-1}{3h}\partial_x(h^3\partial_x \Psi)+ \dfrac{1}{2h}[\partial_x(h^2(\partial_x b) \Psi)-h^2(\partial_x b)(\partial_x \psi)]+(\partial_x b)^2 \Psi .$

At this point, it is worth noting that a rigorous justification (consistency, well-posedness, and stability) of the GN model \eqref{eq1} has been provided in~\cite{DucheneIsrawiTalhouk16} in the flat topography case with a modified velocity variable. The goal of this paper is to consider the more general configuration of variable topography, which has not yet been investigated.

\subsection{The medium amplitude, large topography two-fluid system (GNCH)}\label{section3}
In this section, we build a new (modified) GNCH model ($\epsilon=\OO(\sqrt{\mu})$) from the ``original" model~\eqref{eq1}. Taking into account large topographic variations, we relax the smallness assumption on the amplitude topographic variations parameter $\beta\sim1$, contrary to~\cite{IsrawiLteifTalhouk15,LteifIsrawi18}. This makes more sense in ocean beds. In fact, we assume that there is a $\beta_{max}< \infty$ such that
\begin{equation}
\beta=\OO(1) \hspace{1cm}   \mbox {with} \hspace{1cm}   \beta\in[0,\beta_{max}].\nn
\end{equation} 
Using the following asymptotic expansion: $\frac{1}{1-X}=1+X+\OO(X^2)$ with $X\ll 1$, one gets:
%
\begin{align*}
\dfrac{h_1}{h_1+\gamma h_2}&= \dfrac{1}{\delta(1-\gamma\beta b)+\gamma}\Big(\delta-\delta\epsilon \zeta +\dfrac{\epsilon\zeta\delta^2(1-\gamma)}{\delta(1-\gamma\beta b)+\gamma} + \OO(\epsilon^2)\Big) \;, \\\vspace{1mm}
\dfrac{h_2}{h_1+\gamma h_2} &= \dfrac{1}{\delta(1-\gamma\beta b)+\gamma}\Big(1+\delta\epsilon \zeta-\delta\beta b +\dfrac{(1-\delta\beta b)\epsilon\zeta(\delta-\gamma\delta)}{\delta(1-\gamma\beta b)+\gamma} + \OO(\epsilon^2)\Big)\;.\end{align*}
After replacing these functions with their corresponding approximations in \eqref{qu1} and \eqref{rbar}, the following are obtained:
 \begin{align}\label{Q}
        \overline{\Q}[\epsilon\zeta,\beta b] v  &=  - \lambda(\beta b)\partial_x^2 v\nn\\
        &\quad+\epsilon \left(\theta(\beta b)v \partial_x^2\zeta +  \big(2\theta(\beta b)+(\gamma-1)g(\beta b)\big)\partial_x\zeta\partial_x v+\big(\theta(\beta b)+\dfrac23 (\gamma-1)g(\beta b)\big)\zeta\partial_x^2 v\right)\nn\\ &\quad+\beta \left(\alpha(\beta b)v \partial_x^2 b +  2\alpha(\beta b)\partial_x b\partial_x v+\big(\dfrac\gamma3f(\beta b)+\dfrac23\delta^{-1}f(\beta b)\big)b\partial_x^2 v\right)\nn 
        \\ &\quad+ \epsilon\beta\Big(\big(\theta_1(\beta b)-\alpha_1(\beta b)\big)\zeta v \partial_x^2 b+\big(2\theta_1(\beta b)-\alpha_1(\beta b)\big)\zeta\partial_x b\partial_x v\Big)\nn\\ & \quad+ \epsilon\beta\left(\big(2\theta_1(\beta b)-2\alpha_1(\beta b)+\dfrac13(\delta^{-1}-\beta b)f'(\beta b)-\dfrac\gamma3 g'(\beta b)\big)\partial_x\zeta\partial_x b v\right)\nn
        \\ &\quad+\beta^2\left(\eta(\beta b)(\partial_x b)^2 v +\dfrac{2\gamma}{3}f'(\beta b)b\partial_x b 
        \partial_x v +\dfrac\gamma3 f'(\beta b)b\partial_x^2 b v -\dfrac13f(\beta b)b^2 \partial_x^2 v\right)\nn
        \\ & \quad+\epsilon\beta^2\left(\eta_1(\beta b)(\partial_x b)^2\zeta v\right)
        +\beta^3\left(\dfrac\gamma3 f''(\beta b) (\partial_x b)^2 b v \right)
         +  \OO(\epsilon^2) \; ,\end{align}        
\begin{equation} \label{R}       \overline{\R}[\epsilon\zeta,\beta b] v =  (1-\gamma)g(\beta b)^2\left(\frac12 (\partial_x  v)^2+\frac13 v\partial_x^2 v\right)+ s(\beta b) v^2 + t(\beta b) v\partial_x v + \OO(\epsilon)\;.\qquad \qquad \qquad \quad \;
 \end{equation}
The purely topographical functions that appear in the two preceding expressions are as follows:
\begin{align*}
{\lambda(\beta b)}&=\frac{1+\gamma\delta}{3\delta(\delta(1-\gamma\beta b)+\gamma)}, \qquad{f(\beta b)}=\frac{\delta}{\delta(1-\gamma\beta b)+\gamma},\qquad{g(\beta b)}=\frac{1-\delta\beta b}{\delta(1-\gamma\beta b)+\gamma} \;,\\
 \theta(\beta b) & =\dfrac13(\delta^{-1}-\beta b)^2 f(\beta b) -\dfrac13(\delta^{-1}-\beta b)^2(1-\gamma)f(\beta b)^2-\dfrac\gamma3f(\beta b)-\dfrac\gamma3 f(\beta b)g(\beta b)(1-\gamma)\;,\\
  \alpha(\beta b)&=-\dfrac13(\delta^{-1}-\beta b)^2f'(\beta b)+\dfrac{(\delta^{-1}-\beta b)}{2}f(\beta b)-\dfrac\gamma3\delta^{-1}f'(\beta b)+\dfrac\gamma3 f(\beta b)\;,\\
\theta_1(\beta b)&=\dfrac13(\delta^{-1}-\beta b)^2 f'(\beta b) -\dfrac13(\delta^{-1}-\beta b)^2(1-\gamma)2f(\beta b)f'(\beta b)-\dfrac{(\delta^{-1}-\beta b)}{2}f(\beta b)\nn\\
&\quad + \dfrac{(\delta^{-1}-\beta b)}{2})(1-\gamma)f(\beta b)^2+\dfrac{f(\beta b)}{2}-\dfrac\gamma3 f'(\beta b)+\dfrac{2\gamma}{3} g'(\beta b)-\dfrac23(\delta^{-1}-\beta b)f'(\beta b)\;,\\
\alpha_1(\beta b)&=\dfrac\gamma3(1-\gamma)\Big(g(\beta b)f'(\beta b)+g'(\beta b)f(\beta b)\Big) \;,\\
 \eta(\beta b)&=-\dfrac13(\delta^{-1}-\beta b)^2f''(\beta b)+(\delta^{-1}-\beta b)f'(\beta b)-\dfrac\gamma3\delta^{-1}f''(\beta b)+\dfrac{2\gamma}{3}f'(\beta b) \; , \\
 \eta_1(\beta b)&= \dfrac13(\delta^{-1}-\beta b)^2f''(\beta b)-(\delta^{-1}-\beta b)f'(\beta b)+2(\delta^{-1}-\beta b)(1-\gamma)f'(\beta b)f(\beta b)\nn\\
 &\quad-\dfrac13(\delta^{-1}-\beta b)^2(1-\gamma)4(f'(\beta b))^2-\dfrac13(\delta^{-1}-\beta b)^2(1-\gamma)2f(\beta b)f''(\beta b)\nn\\
 &\quad+ \dfrac{(\delta^{-1}-\beta b)}{3}f''(\beta b)+2f'(\beta b)-(\delta^{-1}-\beta b)f''(\beta b)-\dfrac\gamma3 f''(\beta b)-\dfrac{\gamma(1-\gamma)}{3}f''(\beta b)g(\beta b)\nn\\
 &\quad- \dfrac{\gamma(1-\gamma)}{3}2f'(\beta b)g'(\beta b)-\dfrac{\gamma(1-\gamma)}{3}f(\beta b)g''(\beta b)+\dfrac{2\gamma}{3}g^{(2)}(\beta b)-f'(\beta b) \; ,\\
  s(\beta b)&=\dfrac12(\delta^{-1}-\beta b)^2(\partial_x f(\beta b))^2-2(\delta^{-1}-\beta b)\partial_x(f(\beta b))f(\beta b) \beta \partial_x b\nn \\
 & \quad+\dfrac{(\delta^{-1}-\beta b)^2}{3}f(\beta b)\partial_x^2 (f(\beta b))+\dfrac{\beta^2}{2}(\partial_x b)^2 f(\beta b)^2-\dfrac{\beta}{2}(\delta^{-1}-\beta b)\partial_x^2 b f(\beta b)^2\nn \\
 &\quad-\dfrac\gamma3g(\beta b)\partial_x^2(g(\beta b))-\dfrac\gamma2(\partial_x(g(\beta b)))^2 \; ,\\
t(\beta b)&=\frac53(\delta^{-1}-\beta b)^2 f(\beta b)\partial_x (f(\beta b))-2(\delta^{-1}-\beta b)f(\beta b)^2 \beta \partial_x b-\dfrac{5\gamma}3 g(\beta b)\partial_x (g(\beta b)) \;.
\end{align*}
\begin{remark}
For the sake of readability, we omit the above functions' dependence on $\beta b$ for the rest of the paper. It is also worth noting that the prime notations here and in what follows represent derivatives with respect to $\beta b$, for instance,  $f'(\beta b)=\dfrac{\partial f}{\partial (\beta b)} $. In fact, using the chain rule, one has $\dfrac{\partial}{\partial_x} f(\beta b)=  f'(\beta b)\beta \partial_x b$.\\
We will see later on that the assumption of non-zero depths $(h_1, h_2 >0)$ ensures that the above functions depending on bottom topography are well defined. In fact, $h_1+h_2 >0$ corresponds to $1+1/\delta-\beta b>0$. Since $\gamma<1$, one has $ \delta(1-\gamma\beta b)+\gamma>0$.
\end{remark}

The key to building a simplified asymptotic Green-Naghdi model is plugging the expansions of $\overline{\Q}[\epsilon\zeta,\beta b] v$ and $\overline{\R}[\epsilon\zeta,\beta b] v$ (\eqref{Q}-\eqref{R}) into system~\eqref{eq1}. In fact, all terms of order $\OO(\mu\epsilon^2)$ are now ignored. However, due to the large topography variation assumption, additional topographic terms appear in the expansion of $\overline{\Q}[\epsilon\zeta,\beta b] v$ and $\overline{\R}[\epsilon\zeta,\beta b] v$. Actually, these terms are accompanied by derivative terms on the interface deformation function $\zeta$ and are not controlled in the functional space we are considering. To address this issue, we will construct an equivalent model that can deal with these terms while having a hyperbolic quasilinear structure that allows for full justification.

\subsection{The modified system to be studied}\label{sec3.3}
We start by introducing a specific symmetric differential operator:
      \begin{equation}\label{defT}
      {\mathfrak T}[\epsilon\zeta,\beta b]\Psi \ = \  q_1(\epsilon\zeta,\beta b)\Psi \ + \ \mu\epsilon\beta\kappa_0 \partial_x \zeta \partial_x b \Psi \ - \ \mu\partial_x\Big(\nu q_2(\epsilon\zeta,\beta b)\partial_x\Psi \Big),
      \end{equation}
with $q_1(X,Y)\equiv 1+\kappa_1 X +(\omega_1+\epsilon \phi_1) Y $ and $q_2(X,Y)\equiv 1+\kappa_2 X +\omega_2 Y$ where $\nu,\kappa_0,\kappa_1,\kappa_2,\omega_1, \omega_2$ and $\phi_1$ are functions  of $\beta b$ to be determined in an appropriate way treating all terms that cannot be controlled by the intended energy norm. For the sake of readability, we omit here and for the rest of the paper the dependence of these functions on $\beta b$, thus one can write:
  \begin{align}\label{TversusQ}
        & {\mathfrak T}[\epsilon\zeta,\beta b]\partial_t  v \ -\  q_1(\epsilon\zeta,\beta b) \partial_{ t}\Big(  v \ + \ \mu\overline{\Q}[\epsilon\zeta,\beta b] v \Big)+ q_1(\epsilon\zeta,\beta b) \mu \frac{\gamma+\delta}{\bo}\partial_x^3\zeta \nn \\  &=\mu\epsilon\beta \kappa_0 \partial_x \zeta \partial_x b \partial_t v-\mu\partial_x\nu\partial_x\partial_t v-\mu\nu\partial_x^2\partial_t v-\mu\partial_x(\nu\epsilon\kappa_2 \zeta \partial_x\partial_t v)-\mu\partial_x(\nu\beta\omega_2 b \partial_x\partial_t v) \nn \\
        &\quad -\mu\partial_t \overline{\Q}[\epsilon\zeta,\beta b] v -\mu\epsilon\kappa_1 \zeta \partial_t \overline{\Q}[\epsilon\zeta,\beta b] v-\mu\beta(\omega_1+\epsilon\phi_1) b \partial_t \overline{\Q}[\epsilon\zeta,\beta b] v \nn\\ 
        & \quad+\mu \frac{\gamma+\delta}{\bo}\partial_x^3\zeta +\mu\epsilon\kappa_1\frac{\gamma+\delta}{\bo}\zeta\partial_x^3\zeta+\mu\beta(\omega_1+\epsilon \phi_1) b \frac{\gamma+\delta}{\bo}\partial_x^3\zeta.
      \end{align}
Now we'll look at how to deal with uncontrolled terms that arise in \eqref{TversusQ}. In other words, terms involving third order derivatives on $\zeta$ or $v$, second order derivatives on $\zeta$ and $(\partial_x\zeta)^2$ are uncontrolled.\\
The first order ($\mu\partial_x^2\partial_tv$) terms can be treated by canceling them out using a suitable function $\nu$. In fact, the system's momentum equation~\eqref{eq1} gives 
      \[ \partial_t   v=-(\gamma+\delta)\partial_x\zeta  \ - \ \frac{\epsilon}{2} \partial_x\Big((f^2-\gamma g^2)| v|^2\Big) +\OO(\epsilon^2,\mu).\]
Equivalently, one has
      \begin{equation}\label{zetaxxx} \mu\frac{\gamma+\delta}{\bo}\partial_x^3\zeta \ = \  \frac{-\mu}{\bo}\partial_x^2\partial_t v- \frac{\mu\epsilon}{2\bo} \partial_x^3\Big((f^2-\gamma g^2)| v|^2\Big) +\OO(\mu\epsilon^2,\mu^2).\end{equation}
In equation~\eqref{TversusQ}, substitute the term $\mu\frac{\gamma+\delta}{\bo}\partial_x^3\zeta$ by \eqref{zetaxxx} and applying the time partial derivative to the expansion of $\overline{\Q}[\epsilon\zeta,\beta b] v$ in~\eqref{Q}, 
      thus a suitable choice of $\nu$ is as follows
         \begin{equation}\label{defnu}
 \nu=\lambda-\frac1{\bo}.
      \end{equation}
Using the two approximations that result from~\eqref{eq1}: $$ \partial_t  v=-(\gamma+\delta)\partial_x\zeta  \ +\OO(\epsilon,\mu) \quad \text{and} \quad  \partial_t\zeta=-\partial_x (g v)+\OO(\epsilon,\mu),$$ we will see below that all terms involving $\partial_x^2 \zeta$, $\partial_x^3 \zeta$ and  $(\partial_x \zeta)^2$ can be canceled with an appropriate choice of the functions $\kappa_0, \kappa_1, \kappa_2, \omega_1$, $\omega_2$ and $\phi_1$.\\
For canceling the $b \partial_x^3 \zeta$ terms, the suitable choice of function $\phi_1$ is as follows:
  \begin{equation}\label{defphi1}
 \phi_1 \ = \ \dfrac{\nu\omega_2+\frac\gamma3f+\frac23\delta^{-1}f-\frac{\beta b}{3}f}{ \epsilon I(\beta b)} -\dfrac{\omega_1}{\epsilon} .  \end{equation}
  where $I(\beta b)=\nu-\beta b\big(\frac\gamma3 f+ \frac23\delta^{-1}f \big)+\frac{\beta^2b^2}{3}f$. 
For canceling the $\partial_x b \partial_x^2 \zeta$ terms, we determine the function $\omega_2$ as a solution of the following first order linear differential equation:
\begin{equation}\label{defomega2}
(\nu+\beta b \nu ') \omega_2+\beta b \nu \omega'_2 \ = \ -\nu'-\big(1+\beta (\omega_1+\epsilon \phi_1) b\big)\big(2\alpha +\frac{2\gamma}{3}\beta b f'\big).
  \end{equation}
Furthermore, as $\omega_1 + \epsilon \phi_1$ can be expressed in terms of $\omega_2$ from \eqref{defphi1}, the differential equation~\eqref{defomega2} can be easily rewritten as a first order linear differential equation for $\omega_2$ as follows:
 \begin{eqnarray}\label{ode1}
\omega_2' \ + \  \left(\frac1{\beta b}+\frac{\nu'}{\nu} +\dfrac{ (2 \alpha +\beta b \frac{2\gamma}{3} f')}{I(\beta b)}\right)\omega_2 = -\dfrac{\nu'}{\beta b \nu} -\dfrac{2 \alpha +\beta b \frac{2\gamma}{3} f'}{\beta b \ I(\beta b)}.  \end{eqnarray}
 The first order linear differential equation~\eqref{ode1} can be solved under certain restrictions, to be specified later, on the bottom deformation function $b(x)$, see~\eqref{condbot} in Remark~\ref{oderem}.\\
For canceling the $\epsilon(\partial_x \zeta)^2 \partial_x b$ terms,  $\kappa_0$ is determined as follows:
      \begin{equation}\label{defkappa0}
 \kappa_0 \ = \   \big(2\theta_1-2\alpha_1+\frac13(\delta^{-1}-\beta b)f'-\frac\gamma3 g'\big) \big(1+\beta b\omega_1\big)   \; .
  \end{equation}
  For canceling the $\epsilon\partial_x \zeta \partial_x^2 \zeta$ terms,  $\kappa_2$ is determined as follows:
  \begin{equation}\label{defkappa2}
 \kappa_2 =-\dfrac{\big(3\theta+(\gamma-1)g\big)\big(1+\beta b\omega_1 \big)}{\nu}   \; .
  \end{equation}
    For canceling the $\epsilon\zeta\partial_x^3 \zeta$ terms, with \eqref{defkappa2} in hands, $\kappa_1$ is determined as follows:
\begin{equation}\label{defkappa1}
 \kappa_1 \ = \ \dfrac{\big(-2\theta-\frac{(\gamma-1)}{3}g\big)\big(1+\beta b\omega_1  \big)}{I(\beta b)}   \;.
  \end{equation}
For canceling the $\epsilon\zeta\partial_x b \partial_x^2 \zeta$ terms, we set the function $\omega_1$ so that it satisfies the below equation:
\begin{equation}\label{defomega1}
(\nu\kappa_2)'+(2\theta_1-\alpha_1)(1+\beta b\omega_1)+ \kappa_1\big(\beta b \frac{2\gamma}{3}f'+ 2 \alpha\big)=0.\end{equation}
Fuethermore using the definitions of $\kappa_2$ and $\kappa_1$ in~\eqref{defkappa2} and~\eqref{defkappa1} respectively and after straightforward computations, the above equation turns out to be a first order linear differential equation for $\omega_1$ as follows:
 \begin{align}\label{ode2}
\omega_1' \ + \  \left(\frac{3\theta'+(\gamma-1)g'}{3\theta+(\gamma-1)g}+\frac{1}{\beta b}-\dfrac{2\theta_1-\alpha_1}{3\theta+(\gamma-1)g}+\dfrac{\big(2\theta+\frac{(\gamma-1)}{3}g\big)\big( 2 \alpha +\beta b \frac{2\gamma}{3} f'\big)}{(3\theta+(\gamma-1)g)I(\beta b)}\right)\omega_1\nn \\ = -\frac{3\theta'+(\gamma-1)g'}{\beta b(3\theta+(\gamma-1)g)}+\dfrac{2\theta_1-\alpha_1}{\beta b(3\theta+(\gamma-1)g)}-\dfrac{\big(2\theta+\frac{(\gamma-1)}{3}g\big)\big( 2 \alpha +\beta b \frac{2\gamma}{3} f'\big)}{\beta b(3\theta+(\gamma-1)g)I(\beta b)}.\end{align}
Once more, the first order linear differential equation~\eqref{ode2} can be solved under certain restrictions, to be specified later, on the bottom deformation function $b(x)$, see~\eqref{condbot} in Remark~\ref{oderem}.

  \begin{remark}\label{oderem}
  A first order linear differential equation has the following form: \begin{equation*}y'+n(x) y =p(x),\quad \mbox{where} \ \ y=y(x)\quad \mbox{and} \quad y'=\dfrac{dy}{dx},\end{equation*} where $n(x)$ and $p(x)$ must be continuous functions. The general solution is given by: $y= C e^{-F(x)}+ e^{-F(x)} \int e^{F(x)} p(x) dx,  \ \mbox{where} \ F(x)=\int n(x) dx$, and $C$ is an arbitrary constant. In order to solve~\eqref{ode1} and~\eqref{ode2}, the continuity of the corresponding functions $n(x)$ and $p(x)$ is sufficient. Moreover, for $\kappa_1$, $\kappa_2$, $\omega_1$, $\omega_2$ and $\phi_1$ to be well-defined, altogether, this require certain limitation conditions consisting of additional assumptions on the bottom deformation function $b(x)$. Let us briefly present these conditions.
As a matter of fact, one needs $\beta b$, $\nu$, $I(\beta b)$ and $3\theta+(\gamma-1)g$ to be all different from zero. Since $\gamma\geq 0$, $\delta>0$ and $\bo\gg 6$ then the discriminant of $I(\beta b)$ is $\Delta=\big[\gamma^2\delta^4\bo(\bo-6)+12\delta^4\bo+9\gamma^2\delta^4\big]>0$ thus one needs:
\begin{align}\label{condbot}\beta b \neq 0, \quad \nn \nu(\beta b )\neq 0, \quad \nn  \beta b \neq \displaystyle{\frac{2\delta\bo+\gamma\delta^2\bo-3\gamma\delta^2}{2\delta^2\bo}\pm\frac{ \sqrt{\Delta}}{2\delta^2\bo}}, \quad  \nn \beta b \neq \dfrac{1}{\delta}+\gamma\nn  \quad  \text{and} \quad \beta b \neq \dfrac{1}{\delta}\pm \dfrac{1}{\sqrt{\gamma}} \nn  \tag{H0}.\end{align}
The third and fourth conditions ensure that $I(\beta b)\neq 0$ and the last three conditions ensure that $3\theta+(\gamma-1)g \neq 0$. The above conditions are sufficient to ensure the continuity of the corresponding functions $n(x)$ and $p(x)$ in both~\eqref{ode1} and~\eqref{ode2}.
 \end{remark}
 
Now, handling the $\partial_x^3 v$ terms that remain in~\eqref{TversusQ} and in $\partial_x\big(\overline{\R}[\epsilon\zeta,\beta b]  v\big)$ that are not controlled by the intended energy norm is done by introducing a new function of $\beta b$ to be precisely determined, denoted by $\varsigma$ and embedded in the term $ {\mathfrak T}[\epsilon\zeta,\beta b](\epsilon\varsigma   v\partial_x  v)$. In fact, one has
      \begin{align}\label{TversusR}
        {\mathfrak T}[\epsilon\zeta,\beta b](\epsilon \varsigma  v\partial_x   v) \ +\ \mu\epsilon q_1(\epsilon\zeta,\beta b) \partial_{x}\Big(\overline{\R}[\epsilon\zeta,\beta b]  v\Big)\nn\\
       \quad = \  q_1(\epsilon\zeta,\beta b)(\epsilon \varsigma  v\partial_x  v) \ + \ \mu\epsilon\beta \kappa_0 \partial_x \zeta \partial_x b\big(\epsilon\varsigma v\partial_x v)\nn\\ \ - \ \mu \partial_x \Big(\nu q_2(\epsilon\zeta,\beta b)\partial_x (\epsilon\varsigma v\partial_x  v ) \Big)+\mu\epsilon q_1(\epsilon\zeta,\beta b )\partial_x\Big(\overline{\R}[\epsilon\zeta,\beta b]  v\Big).
      \end{align}
Adding~\eqref{TversusQ} to~\eqref{TversusR}, and fixing $\varsigma$ as follows 
   \begin{align}\label{defvarsigma}
	  \nu(1+\beta b \omega_2)\varsigma  &=  \frac{(1-\gamma)g^2}{3}(1+\beta\omega_1 b)-\frac{1}{\bo}(f^2-\gamma g^2)+\theta g(1+\beta\omega_1 b)\nn\\ &\quad +\nu\beta b \omega_2(f^2-\gamma g^2)+(1+\beta\omega_1 b)\big[\beta b (\frac\gamma 3f+\frac23\delta^{-1} f)(f^2-\gamma g^2)\big]\nn\\ &\quad -(1+\beta\omega_1 b)\big[\beta^2 b^2\frac13f(f^2-\gamma g^2)\big]-\beta w_1 b \lambda(f^2-\gamma g^2) ,
      \end{align}
     while dropping all terms of order $  \OO(\mu^2,\epsilon^2\mu)$ yields the following approximation:
\begin{align}\label{eqcons}
       {\mathfrak T}(\partial_t  v +\epsilon \varsigma  v\partial_x  v) - q_1(\epsilon\zeta,\beta b) \partial_{ t}\Big(  v \ + \ \mu\big(\overline{\Q}[\epsilon \zeta, \beta b] v\big) \Big)\nn \\+q_1(\epsilon\zeta,\beta b) \mu \frac{\gamma+\delta}{\bo}\partial_x^3\zeta+ \mu\epsilon q_1(\epsilon\zeta,\beta b) \partial_{ x}\Big(\overline{\R}[\epsilon \zeta, \beta b]  v   \Big) \\
         = q_1(\epsilon\zeta,\beta b) (\epsilon \varsigma v \partial_x v) +\mu[\mathcal{A}]v\partial_x v + \mu [\mathcal{B}](\partial_x v)^2 + \mu[\mathcal{C}]v \partial_x^2 v\nn\\+\mu[\mathcal{D}]\partial_x\big((\partial_x v)^2\big)+\mu[\mathcal{E}]v^2+\mu[\mathcal{F}]\partial_x \zeta+\OO(\mu^2,\mu\epsilon^2)\nn. \end{align}
We would like to point out that the terms $\mathcal{A}$, $\mathcal{B}$, $\mathcal{C}$, $\mathcal{D}$, and $\mathcal{E}$ only represent functions that depend on $b(x)$, whereas the term $\mathcal{F}$ represents a function that depends on both $\zeta(t,x)$ and $b(x)$. We do not try to give in here their exact expressions for the sake of simplicity. These functions, however, are described in detail at the end of this paper in Appendix~\ref{appendixB}. In fact, their precise expressions are irrelevant to our current purpose. Despite their length, these functions are easy to control and are accompanied by controllable terms.

Now, the last step to get the new equivalent asymptotic model is to multiply the second equation of~\eqref{eq1} by $ q_1(\epsilon\zeta,\beta b)$ and include the approximation~\eqref{eqcons}. This yields the following GNCH model with large topography:
      \begin{equation}\label{eq:Serre2}\left\{ \begin{array}{l}
      \partial_{ t}\zeta +\partial_x\left(H(\epsilon\zeta,\beta b)v\right)\ =\  0, \vspace{2mm}\\ 
      \mathfrak T[\epsilon\zeta,\beta b] \left( \partial_{ t}   v + \epsilon\varsigma { v } \partial_x {  v } \right) + (\gamma+\delta)q_1(\epsilon\zeta,\beta b)\partial_x
      \zeta +\frac\epsilon2 q_1(\epsilon\zeta,\beta b) \partial_x \left(H'(\epsilon\zeta,\beta b)| v|^2\right)-q_1(\epsilon\zeta,\beta b) (\epsilon \varsigma v \partial_x v)\vspace{1mm} \\
      = \mu [\mathcal{A}]v\partial_x v + \mu [\mathcal{B}](\partial_x v)^2 + \mu[\mathcal{C}]v \partial_x^2 v +\mu[\mathcal{D}]\partial_x\big((\partial_x v)^2\big)+\mu[\mathcal{E}]v^2+\mu[\mathcal{F}]\partial_x \zeta \;,
      \end{array} \right. \end{equation}
where $H(\epsilon\zeta,\beta b)  = \dsp\frac{h_1h_2}{h_1+\gamma h_2}$, 
$H'(\epsilon\zeta,\beta b)  = \dsp\frac{h_1^2-\gamma h_2^2}{(h_1+\gamma h_2)^2}$, $h_1=1-\epsilon\zeta$, $h_2=1/\delta+\epsilon\zeta-\beta b$, $q_1(\epsilon\zeta,\beta b)\equiv 1+\kappa_1 \epsilon \zeta +(\omega_1+\epsilon \phi_1) \beta b $ and $q_2(\epsilon\zeta,\beta b)\equiv 1+\kappa_2 \epsilon \zeta +\omega_2 \beta b$
with $\nu,\kappa_0,\kappa_1,\kappa_2,\omega_1,\omega_2$ and $\phi_1$ are functions of $\beta b$ defined in~\eqref{defnu}~\eqref{defkappa0}~\eqref{defkappa1}~\eqref{defkappa2},~\eqref{defomega1},~\eqref{defomega2},~\eqref{defphi1}.

\begin{remark}
It is worth mentioning that the same model derived in~\cite{DucheneIsrawiTalhouk15} can be easily recovered by setting $\beta=0$ in~\eqref{eq:Serre2}. Moreover, the same model derived in~\cite{IsrawiLteifTalhouk15}  can be recovered by setting $\beta=\OO(\sqrt{\mu})$ and dropping all terms of order $\OO(\mu^2,\mu\epsilon^2,\mu\epsilon\beta,\mu\beta^2)$ in~\eqref{eq:Serre2}. Finally,  the same model derived in \cite{LteifIsrawi18} can be recovered by considering the bottom deformation function $b=b^{(\alpha)}(x)=b(\alpha x)$ and assuming $\beta \alpha = \OO(\sqrt{\mu})$ while dropping all terms of order $\OO(\mu^2,\mu\epsilon^2,\mu\epsilon\beta\alpha,\mu\beta^2\alpha^2)$. 
\end{remark}
       
\section{Properties of operator $\mathfrak{T}$ and its inverse}\label{Sec4}

Let us now present the following assumptions, which are critical for the mathematical analysis:\\ 
$\bullet$ Assume that $\nu(\beta b)>0$ has a lower bound: \begin{equation*}\nu(\beta b)=\lambda(\beta b)-\dfrac{1}{\bo}\geq\nu_{0}>0.\end{equation*}
Hereinafter, this condition is added to the previous assumptions on the bottom deformation function $b(x)$, namely~\eqref{condbot}.\\
$\bullet$ Assume that there exist positive constant $h_{min}>0$ such that
\begin{equation}\label{CondDepth}\tag{H1}
\min\Big\lbrace  \inf_{x\in \RR} h_1= \inf_{x\in \RR} [1-\epsilon\zeta],\inf_{x\in \RR} h_2= \inf_{x\in \RR} [1/\delta -\epsilon\zeta +\beta b] \Big\rbrace\geq h_{min} \; .
\end{equation}
In the oceanographic context, the previous assumptions are referred to as the depth-condition on both layers of the fluid.\\
$\bullet$ Assume that $q_1(\epsilon\zeta,\beta b)+\mu\epsilon\beta\kappa_0 \partial_x \zeta \partial_x b >0$ and $q_2(\epsilon\zeta,\beta b)>0$ have lower bounds:
\begin{equation}\label{CondEllipticity}\tag{H2}
\exists\ q_{min}>0 , \mbox{ such that } \quad \min\Big\lbrace \inf_{x\in \RR}\Big(q_1(\epsilon\zeta,\beta b)+\mu\epsilon\beta\kappa_0 \partial_x \zeta \partial_x b\Big)  , \inf_{x\in \RR}  q_2(\epsilon\zeta,\beta b) \Big\rbrace \geq q_{min} \;. 
\end{equation}
As a matter of fact, for fixed $\vert \zeta \vert _{W^{1,\infty}}$ and $\vert b\vert _{W^{1,\infty}}$, the two conditions~\eqref{CondDepth} and~\eqref{CondEllipticity} reduce to an estimate on
\begin{align*}
&\max(|\kappa_1|_{L^\infty},|\kappa_2|_{L^\infty},1,\delta_{\max})\epsilon_{\max}\big\vert \zeta\big\vert_{L^\infty}\\
+&\max(|\omega_1|_{L^\infty},\epsilon|\phi_1|_{L^\infty},|\omega_2|_{L^\infty}, \delta_{max})\beta_{max}\big\vert b\big\vert _{L^\infty} +\mu_{\max}\epsilon_{\max}\beta_{\max}\Big(\big\vert \kappa_0\big\vert_{L^\infty}\big\vert \partial_x \zeta \big\vert_{L^\infty}\vert \partial_x b\big\vert_{L^\infty}) \le 1-h_{min},
\end{align*}
with $\epsilon_{\max}=\min(M\sqrt{\mu_{\max}},1)$. In that case,~\eqref{CondDepth} and~\eqref{CondEllipticity} hold for any parameter in $ \P_{\rm CH}$.

This section will expose the strict dependence of the assumption~\eqref{CondEllipticity} to our framework. In fact, much of the analysis hinges on the left-most operator $\mfT$ introduced in~\eqref{defT} that we shall recall 
\begin{equation}\label{def22}
{\mathfrak T}[\epsilon\zeta,\beta b]V \ = \ q_1(\epsilon\zeta,\beta b)V \ + \ \mu\epsilon\beta\kappa_0 \partial_x \zeta \partial_x b  V \ - \ \mu\partial_x\Big(\nu q_2(\epsilon\zeta,\beta b)\partial_xV \Big),
\end{equation}
where $q_1(\epsilon\zeta,\beta b)\equiv 1+\kappa_1 \epsilon \zeta +(\omega_1+\epsilon \phi_1) \beta b $ and $q_2(\epsilon\zeta,\beta b)\equiv 1+\kappa_2 \epsilon \zeta +\omega_2 \beta b$
with $\nu,\kappa_0,\kappa_1,\kappa_2,\omega_1,\omega_2$ and $\phi_1$ are functions of $\beta b$ defined in~\eqref{defnu}~\eqref{defkappa0}~\eqref{defkappa1}~\eqref{defkappa2},~\eqref{defomega1},~\eqref{defomega2},~\eqref{defphi1}.  
As we will see later in section \ref{The-linearized-system}, in order to write system \eqref{eq:Serre2} as a quasilinear hyperbolic system, the inverse operator, namely $\mfT^{-1}$, must be applied to the latter system's second equation. Thus, one must first ensure that the strong ellipticity property of $\mfT$ is maintained so that its inverse $\mfT^{-1}$  is well-posed and continuous.

 As in  \cite{LteifIsrawi18,IsrawiLteifTalhouk15}, we will work in the space $H^1_\mu(\RR)$ endowed with the norm \[\forall \  v\in H^1(\RR),\ \quad
\vert\ v\ \vert^2_{H^{1}_\mu}\ =\ \vert\ v\ \vert^2_{L^2}\ +\ \mu\ \vert\ \partial_x v\ \vert^2_{L^2},
\] 
where  $\vert \cdot \vert_{H^{1}_\mu}$ is equivalent to  $\vert \cdot \vert_{H^{1}}$ but not uniformly with respect to $\mu$. The below gives the intvertibility of $\mfT$.
\begin{lemma}\label{proprim}
Assume that  $(\zeta,b)\in W^{1,\infty}(\RR)^2$ under the assumption~\eqref{CondEllipticity}.
Then the operator
\[
{\mfT}: H^2(\RR)\longrightarrow L^2(\RR)
\]
is one-to-one and onto.
\end{lemma}
\begin{proof}
The operator's invertibility is an application of the Lax-Milgram theorem, as stated in~\cite[Lemma 4.2]{LteifIsrawi18}. Using~\eqref{CondEllipticity}, the bilinear form, consisting of:
\begin{align*}
a(u,v)=\big({\mfT} u , v\big)=\big(\ (1+\epsilon \kappa_1\zeta+\beta (\omega_1+\epsilon \phi_1) b +\mu\epsilon\beta \kappa_0 \partial_x \zeta \partial_x b)u\ ,\ v\ \big)+\mu \big(\ \nu(1+\epsilon\kappa_2\zeta+\beta \omega_2 b)\partial_x u\ ,\ \partial_x v\ \big)\;,
\end{align*}
can be easily demonstrated to be uniformly continuous and coercive on $H^1_\mu(\RR) \times H^1_\mu(\RR)$.
Hence, by Lax-Milgram theorem, for every $f \in L^2(\RR)$,
there exists a unique $u\in H^1_\mu(\RR)$
such that, for all $v\in H^1_\mu(\RR)$ we have
$\colon
a(u,v)=\big({\mfT} u , v \big)=(f,v)
$. It remains to prove that $v\in L^2(\RR)$ as in~\cite[Lemma 4.2]{LteifIsrawi18}.
\end{proof}
Before providing some higher-order estimates on $\mfT$ that will be useful in the paper's sequel. We recall the commutator estimate we will use, which was developed by Kato-Ponce \cite{KatoPonce} and recently improved by Lannes \cite{Lannes06}. Specifically, for any $s>3/2$, and $ N \in H^s(\RR),M\in H^{s-1}(\RR)$, one has
\begin{equation}\label{cest}
\big \vert [\Lambda^s, M]N \big \vert_{2} \lesssim \vert M_x\vert_{H^{s-1}}\vert N\vert_{H^{s-1}} \; .
\end{equation}
Furthermore, we will make extensive use of the classical product estimate (see \cite{AlinhacGerard,Lannes06,KatoPonce}). In particular, for any $M,N\in H^s(\RR)$, $s>3/2$, one has
\begin{equation}\label{mest}
\vert MN\vert_{H^s}\lesssim \vert M\vert_{H^s}\vert N\vert_{H^s} \; .
\end{equation}
We will also make extensive use of the continuous embedding $H^s(\RR)\subset W^{1,\infty}(\RR)$ for $s>3/2$.
\begin{lemma}\label{lemma2}
Fix $ s>3/2$ and assumption~\eqref{CondEllipticity}. Suppose that $\zeta \in  H^{s} (\RR)$, $b \in H^{s+2}(\RR) $, and $(u,v)\in  H^{s}(\RR)\times H^{1}(\RR)$. Then it holds that:
\begin{equation}
\label{eq:estComT}
 \big( [\Lambda^s, \mfT[\epsilon\zeta,\beta b] ]u, v\big)   \leq \  C_0 \max(\epsilon,\beta) \vert b \vert_{H^{s+2}}\vert \zeta\vert_{H^s}\vert u\vert_{H^s}  \vert v \vert_{H^{1}_\mu}  \; , 
\end{equation}
\begin{equation}
\label{eq:estLambdaT}
 \big(\Lambda^s \mfT[\epsilon\zeta,\beta b]  u, v \big)  \leq  C_0 \big(  \max(\epsilon,\beta) \vert b \vert_{H^{s+2}} + \vert b\vert_{H^{s+1}}\big)\vert \zeta\vert_{H^s}\vert u\vert_{H^s}  \vert v \vert_{H^{1}_\mu}  \; ,
\end{equation}
where $C_0$ is a constant that depends on $\vert b\vert_{H^s}$ and is derived from the expressions of  $\nu,\kappa_0,\kappa_1,\kappa_2,\omega_1,\omega_2,\phi_1$.
\end{lemma}
\begin{proof}
\underline{Proof of \eqref{eq:estComT}}: Combining \eqref{def22}, the identities 
$
\kappa_0\partial_x \zeta \partial_x b  =\partial_x ( \kappa_0 \zeta \partial_x b) -\partial_x (\kappa_0 )\zeta \partial_x b -\kappa_0  \zeta \partial_x^2 b$, and  $ \partial_x\big([\Lambda^s , f]g\big
)=[\Lambda^s,\partial_x f]g+[\Lambda^s, f]\partial_x g \; 
$, and by integration by parts with \eqref{cest}-\eqref{mest}, it holds that
\begin{align*}
\big(\ [\Lambda^s, {\mfT} ] u\ ,\  v\ \big) & =\epsilon \big(\  [\Lambda^s, \kappa_1 \zeta] u \ ,\  v\ \big)+\beta \big(\  [\Lambda^s, (\omega_1+\epsilon \phi_1) b] u \ ,\  v\ \big) -\mu\epsilon\beta \big(\  [\Lambda^s, \kappa_0 \zeta\partial_x b] u \ ,\ \partial_x v\ \big)\\
&\quad -\mu\epsilon\beta \big(\  [\Lambda^s, \kappa_0 \zeta\partial_x b] \partial_x u \ ,\  v\ \big)   -\mu\epsilon\beta \big(\ [\Lambda^s, \partial_x (\kappa_0) \zeta\partial_x b] u \ ,\ v\ \big)-\mu\epsilon\beta \big( \ [\Lambda^s, \kappa_0\zeta\partial_x^2 b] u \  , \ v \ \big) \\
&\quad +\mu  \big(\ [\Lambda^s,\nu]\partial_x u\ ,\  \partial_x v\ \big)+\mu \epsilon \big(\ [\Lambda^s,\nu\kappa_2\zeta]\partial_x u\ ,\  \partial_x v\ \big) +\mu \beta \big(\ [\Lambda^s,\nu\omega_2 b]\partial_x u ,  \partial_x v\ \big)\\
&\lesssim \max(\epsilon,\beta) \vert b \vert_{H^{s+2}}\vert \zeta\vert_{H^s}\vert u\vert_{H^s}(\vert v \vert_{L^2} +  \sqrt{\mu}\vert v_x \vert_{L^2}) \; .
\end{align*}
\underline{Proof of \eqref{eq:estLambdaT}}: 
By definition, we have
$$
\big( \Lambda^s {\mfT} u ,  v \big)=\big(\ [\Lambda^s, {\mfT} ] u\ ,\  v\ \big)+\big( {\mfT} \Lambda^s u ,  v \big)=I_1+I_2\;.
$$
With \eqref{eq:estComT} in hands, it remains to bound from above $I_2$. Indeed, by integration by parts, it holds that
\begin{align*}
\big(  {\mfT} \Lambda^s u ,  v \big)&=\big( (1+\epsilon \kappa_1\zeta+\beta(\omega_1+\epsilon \phi_1) b) \Lambda^su ,  v \big)
-\mu\epsilon\beta\big( \kappa_0 u \zeta \partial_x b \Lambda^su, \partial_x v \big)
-\mu\epsilon\beta\big( \partial_x(\kappa_0u)\zeta \partial_x b\Lambda^su ,  v \big)\\
&\quad-\mu\epsilon\beta\big( \kappa_0u\zeta \partial_x^2 b \Lambda^su , v \big)
+\mu \big( \nu(1+\epsilon\kappa_2\zeta+\beta\omega_2 b) \Lambda^s\partial_x u,  \partial_x v \big)\\
&\lesssim   \vert b \vert_{H^{s+1}}\vert \zeta\vert_{L^{\infty}}\vert u\vert_{H^s}(\vert v \vert_{L^2} +\sqrt{\mu}\vert v_x \vert_{L^2}) \; .  
\end{align*}
\end{proof}
We conclude this section by introducing a technical estimate on which much of the analysis in the following sections is based, particularly the derivation of the energy estimate. With Lemma \ref{lemma2} and assumption~\eqref{CondEllipticity} in hands,  we do not attempt to provide a proof for the below estimates because it is a direct adaptation of the proofs in~\cite[Lemma 5.3 and Corollary 5.4]{DucheneIsrawiTalhouk15}.
\begin{corollary}\label{col:comwithT}
Fix $ s>3/2$ and assumption~\eqref{CondEllipticity}. Suppose that $\zeta \in  H^{s} (\RR)$, $b \in H^{s+2}(\RR) $, and $(u,v)\in  H^{s-1}(\RR)\times H^{1}(\RR)$. Then, for $m=\{s,s-1\}$, it holds that:
\begin{equation}\label{i-1}
\Vert \mfT^{-1} \psi\Vert_{H^m(\RR)} + \sqrt{\mu} \Vert \partial_x \mfT^{-1}\psi \Vert_{H^m(\RR)} + \sqrt{\mu} \Vert  \mfT^{-1} \partial_x\psi\Vert_{H^m(\RR)}\leq C_1 \vert \psi\vert_{H^m} \; .
\end{equation}
Furthermore, since $\mfT$ is symmetric, the following estimate holds
\begin{equation*}
 \big( \ \big[\Lambda^s, \mfT^{-1}[\epsilon\zeta,\beta b]\big] u\ ,\ \mfT[\epsilon\zeta,\beta b]  v\ \big)   = -\big(\  \big[\Lambda^s,\mfT[\epsilon\zeta,\beta b] \big] {\mfT}^{-1}  u\ ,\  v\ \big)    \le \max(\epsilon,\beta) C_1C_0 \vert \zeta\vert_{H^s}\big\vert u \big\vert_{H^{s-1}} \big\vert v \big\vert_{H^{1}_\mu}
\label{eq:comwithT}\end{equation*}
with $C_1=C(M_{\rm CH},q^{-1}_{min},\big\vert \zeta \big\vert_{H^s},\big\vert b\big\vert_{H^{s+2}})$.
\end{corollary}     
%

\section{Preliminary procedure}\label{The-linearized-system}
To motivate the introduction of the energy norm (see Definition \ref{defispace}), consider the linearized system of \eqref{eq:Serre2} around a reference state $\underline{U}=(\underline{\zeta},\underline{v})^T$:
\begin{equation}\label{SSlsys}
	\left\lbrace
	\begin{array}{l}
	\dsp\partial_t U+A[\underline{U}]\partial_x U+B(\underline{U})=0\vspace{1mm},
        \\
	\dsp U_{\vert_{t=0}}=U_0.
	\end{array}\right.
\end{equation}
This necessitates first writing the system
\begin{equation}\label{eqn:Serre2varf}\left\{ \begin{array}{l}
\dsp \partial_{ t}\zeta +H(\epsilon\zeta,\beta b)\partial_x  v+\epsilon \partial_x\zeta H'(\epsilon\zeta,\beta b)  v -\beta\partial_x b G(\epsilon\zeta,\beta b)v \ =\  0,\\ \\
\dsp {\mathfrak T} \left( \partial_{ t}  v + \frac{\epsilon}{2} \varsigma \partial_x({v }^2) \right) + (\gamma+\delta)q_1(\epsilon\zeta,\beta b)\partial_x \zeta + \epsilon q_1(\epsilon\zeta,\beta b)\big(H'(\epsilon\zeta,\beta b)-\varsigma \big) v\partial_x v\\+ \epsilon q_1(\epsilon\zeta,\beta b)\partial_x \big(\dfrac{H'(\epsilon\zeta,\beta b)}{2}\big) v^2 = \mu [\mathcal{A}]v\partial_x v + \mu [\mathcal{B}](\partial_x v)^2 + \mu[\mathcal{C}]v \partial_x^2 v\\ \hspace{5cm}+\mu[\mathcal{D}]\partial_x\big((\partial_x v)^2\big)+\mu[\mathcal{E}]v^2+\mu[\mathcal{F}]\partial_x \zeta .
\end{array} \right. \end{equation}
such that $G(\epsilon\zeta,\beta b)=\Big(\frac{h_1}{h_1+\gamma h_2}\Big)^2$, $H(\epsilon\zeta,\beta b)  = \dsp\frac{h_1h_2}{h_1+\gamma h_2}$, 
$H'(\epsilon\zeta,\beta b)  = \dsp\frac{h_1^2-\gamma h_2^2}{(h_1+\gamma h_2)^2}$, $h_1=1-\epsilon\zeta$, $h_2=1/\delta+\epsilon\zeta-\beta b$, $q_1(\epsilon\zeta,\beta b)\equiv 1+\kappa_1 \epsilon \zeta +(\omega_1+\epsilon \phi_1) \beta b $ and $q_2(\epsilon\zeta,\beta b)\equiv 1+\kappa_2 \epsilon \zeta +\omega_2 \beta b$
with $\nu,\kappa_0,\kappa_1,\kappa_2,\omega_1,\omega_2$ and $\phi_1$ are functions of $\beta b$ defined in~\eqref{defnu}~\eqref{defkappa0}~\eqref{defkappa1}~\eqref{defkappa2},~\eqref{defomega1},~\eqref{defomega2},~\eqref{defphi1},
 as a quasilinear system of first-order evolution equations, i.e. in the form
\begin{equation}\label{condensedeq}
\partial_tU+A[U]\partial_xU+B(U) = 0 \; .
\end{equation}
Indeed, applying $\mfT^{-1}$ to the second equation of \eqref{eqn:Serre2varf}, the nonlinear condensed equation \eqref{condensedeq} reads:
\begin{equation}\label{defA0A1}
A[U]=\begin{pmatrix}
\epsilon H'(\epsilon \zeta,\beta b) v & H(\epsilon\zeta,\beta b)\\
\mfT^{-1}(Q_0 \cdot+\epsilon^2Q_1\cdot)& \mfT^{-1}(\mfQ[\epsilon\zeta,\beta b,v] \cdot)+\epsilon\varsigma v
\end{pmatrix},
\end{equation}
\begin{equation}\label{defB}
B[U]=\begin{pmatrix}
 -\beta\partial_x b G(\epsilon\zeta,\beta b)v\\
\mfT^{-1}\Big( \gamma\epsilon\beta q_1(\epsilon\zeta,\beta b)h_1(h_1+h_2)(h_1+\gamma h_2)^{-3}\partial_x b v^2-\mu[\mathcal{E}]v^2\Big) 
\end{pmatrix}\;,
\end{equation}
where we denote by
\begin{equation}\label{defQ0Q1}
Q_0\equiv Q_0(\epsilon\zeta,\beta b)  = (\gamma+\delta)q_1(\epsilon\zeta,\beta b)-\mu[\mathcal{F}],\qquad Q_1\equiv Q_1(\epsilon\zeta,\beta b,v)=-\gamma q_1(\epsilon \zeta,\beta b)\dfrac{(h_1+h_2)^2}{(h_1+\gamma h_2)^3}{ v}^2 \;,
\end{equation}
\begin{equation}\label{defmfQ}
\text{and} \quad \mfQ[\epsilon\zeta,\beta b,v]f \ \equiv \ \Big(\epsilon q_1(\epsilon\zeta,\beta b)(H'(\epsilon\zeta,\beta b) -\varsigma)-\mu[\mathcal{A}]\Big)vf-\mu[\mathcal{B}]\partial_xv f -\mu[\mathcal{C}]v\partial_x f-\mu[\mathcal{D}]\partial_x(\partial_x v f).
\end{equation}
We recall that  the expressions of $\beta b$-functions $\mathcal{A}$, $\mathcal{B}$, $\mathcal{C}$, $\mathcal{D}$, $\mathcal{E}$, $\mathcal{F}$ are given in Appendix~\ref{appendixB} with $\mathcal{F}$ in addition depending on $\zeta$. At this stage a good suggestion of a pseudo-symmetrizer for $A[\underline{U}]$ is given by the symmetric matrix:
\begin{equation}\label{defZ}Z[\underline{U}]=\begin{pmatrix}
 \dfrac{\underline{Q}_0+\epsilon^2 \underline{Q}_1}{H(\epsilon\underline{\zeta},\beta b)}& 0 \\
0&\underline{\mfT}
\end{pmatrix},\qquad\text{ where }\qquad  \underline{\mfT}=  q_1(\epsilon\underline{\zeta},\beta b) \ + \ \mu\epsilon\beta\kappa_0 \partial_x \underline{\zeta} \partial_x b   \ - \ \mu\partial_x\Big(\nu q_2(\epsilon\underline{\zeta},\beta b)\partial_x \Big)\;,
\end{equation}
$\underline{Q}_0 = Q_0(\eps\underline{\zeta},\beta b)$, and $\underline{Q}_1= Q_1(\eps\underline{\zeta},\beta b,\underline{v})$.
The assumption below ensure a positive definite pseudo-symmetrizer:
\begin{equation}\label{H3} \tag{H3} \exists \  Q_{min} >0 \  \  \mbox {such that} \  \    \underline{Q}_0+\epsilon^2 \underline{Q}_1\geq Q_{min} > 0.\end{equation}
We normally define the energy of the linearized system \eqref{SSlsys} as:
\begin{equation}\label{es}
 E^s(U)^2=(\Lambda^sU,Z[\underline{U}]\Lambda^sU)=( \Lambda^s\zeta,\frac{\underline{Q}_0+\epsilon^2 \underline{Q}_1}{H(\epsilon\underline{\zeta},\beta b)}\Lambda^s\zeta)+\left(\Lambda^s v, \underline{\mfT}  \Lambda^s v\right) \;.
\end{equation}
	Following the usual assumptions~\eqref{CondDepth},~\eqref{CondEllipticity} and~\eqref{H3}, it is common to apply the following energy space to our problem:
	\begin{definition}[Energy space]\label{defispace}
For any $s\ge0$, $Y^s$ denotes the vector space $H^s(\RR)\times H^{s+1}(\RR)$ endowed with the norm
\begin{equation}\label{endowed-norm}
\vert U\vert^2_{Y^s}\equiv \vert \zeta\vert^2 _{H^s}+\vert v\vert^2 _{H^s}+ \mu\vert \partial_xv\vert^2 _{H^s}\;.
\end{equation}
We refer to $Y^s_T$ the space $C([0, T/\max(\eps,\beta)];Y^{s})$ endowed with its canonical norm.
\end{definition}
Under the additional assumption given in~\eqref{CondDepth},~\eqref{CondEllipticity} and~\eqref{H3}, $E^s(U)$ is equivalent to the $|\cdot|_{Y^s}$. This equivalence is asserted in the following Lemma. We omit the proof of this Lemma because it is demonstrated using the same techniques as in~\cite[Lemma 5.2]{IsrawiLteifTalhouk15}.
\begin{lemma}\label{lemmaes}
Under the additional assumptions~\eqref{CondDepth},~\eqref{CondEllipticity} and ~\eqref{H3} suppose that $(\underline{\zeta},\underline{v})\in W^{1,\infty}(\RR)$ and $b\in H^{s+2}(\RR)$. Then for any $s\ge0$ the norm $\vert \cdot\vert_{Y^s }$ and the natural energy $E^s(U)$ are uniformly equivalent with respect to $\mu,\eps,\beta \in (0,1)^3$, moreover there exist positive constant  $c_0=C(M_{\rm CH},h_{min}^{-1},q_{min}^{-1},Q_{min}^{-1}, \epsilon|\underline{U}|_{W^{1,\infty}},\beta|b|_{W^{2,\infty}})$ such that:
\begin{equation*}
E^s(U) \leq c_0  \vert U \vert_{Y^s} \qquad \text{ and } \qquad \vert U \vert_{Y^s}  \leq c_0 E^s(U) \; .
\end{equation*}
\end{lemma}

\begin{lemma}\label{lemma4}
Suppose that $s\ge0$, $b \in H^{s+3}(\RR)\cap W^{1,\infty}(\RR) $, and $(\underline{\zeta},\underline{v}) \in  H^{s}(\RR)^2\cap W^{1,\infty}(\RR)^2$. Then it holds that:
\begin{align*}
\label{assemb}
&\bullet \big\vert q_1(\epsilon\underline{\zeta},\beta b) -1 \big\vert_{H^s} + \big\vert q_2(\epsilon\underline{\zeta},\beta b) -1 \big\vert_{H^s} \le \max(\eps,\beta) C \big(  \vert b\vert_{H^{s+1}},\vert\underline{\zeta}\vert_{H^s} \big), \\
&\bullet\big\vert \mathcal{A} \big\vert _{H^s}+\big\vert \mathcal{B} \big\vert _{H^s}+\big\vert \mathcal{C} \big\vert _{H^s}+\big\vert \mathcal{D} \big\vert _{H^s}+\big\vert \mathcal{E} \big\vert _{H^s} \le \max(\eps,\beta) C \big(  \vert b\vert_{H^{s+3}}  \big), \\
&\bullet \big\vert \mathcal{F} \big\vert _{H^s} \le \max(\eps,\beta) C \big(  \vert b\vert_{H^{s+2}},\vert\underline{\zeta}\vert_{H^s} \big)\;.
\end{align*}
Moreover, under the assumption \eqref{CondDepth}, it holds that
\begin{align*}
\Big\vert \frac{ \underline{Q}_0\ +\eps^2 \underline{Q}_1}{H(\eps\underline{\zeta},\beta b)} \Big\vert_{L^{\infty} } +  \Big\vert H'(\eps\underline{\zeta},\beta b) \Big\vert_{L^{\infty} }& \le   C \big( h_{min}^{-1}, \vert b\vert_{L^{\infty}},\vert\underline{\zeta}\vert_{L^{\infty}}  \big)  \; ,\\
\Big\vert  \partial_x\big(\dfrac{Q_0(\epsilon\zetabar,\beta b)+\epsilon^2Q_1(\epsilon\zetabar,\beta b,\vbar)}{H(\epsilon\zetabar,\beta b)}H'(\epsilon\zetabar,\beta b)\vbar\big) \Big\vert_{L^{\infty}} &\le C \big( h_{min}^{-1}, \vert b\vert_{W^{1,\infty}(\RR)},\vert\underline{\zeta}\vert_{W^{1,\infty}(\RR)} \big)  \vert\underline{v}\vert_{W^{1,\infty}(\RR)} \; ,
\end{align*}
\end{lemma}

\begin{proof}
It is not hard to check that, by definition and using intensively the product estimate \eqref{mest}, the desired estimates holds.
\end{proof}

\section{Long-term well-posedness}\label{Sec6}
The first part focuses on the mathematical analysis of the linearized system \eqref{SSlsys}. Finally, in the second part, we deduce the main result of this section, which is the well-posedness of the nonlinear system~\eqref{condensedeq} on time scales of order $\max(\eps,\beta)^{-1}$.

\subsection{The linear problem}\label{linear-analysis}
The following subsection targets the existence and uniqueness of solution to the initial value problem \eqref{SSlsys} and thus we focus our attention on the main point of interest, i.e. the $Y^{s>3/2}$ energy estimate \eqref{energy-est-Mex}, which gives uniform estimates on a time scale of order $\max(\eps,\beta)^{-1}$.
\begin{proposition}\label{prop1}
Fix $s>3/2$ and  \eqref{condbot}. Let $T>0$ and $\underline{U}\in Y^s$ satisfying the assumptions \eqref{CondDepth},~\eqref{CondEllipticity} and~\eqref{H3} such that $ \partial_t \underline{\zeta} , \partial_t \underline{v} \in L^{\infty}(\RR)$. For any initial data $U_0\in Y^s$ satisfying the assumptions~\eqref{CondDepth},~\eqref{CondEllipticity} and~\eqref{H3} there exists a unique solution $U=(\zeta, v)^T$ $\in Y^s_{ T }$ to (\ref{SSlsys}) preserving  \eqref{CondDepth} for all $0\leq \max(\eps,\beta)t\leq T$, such that the energy estimates hold
\begin{equation}\label{energy-est-Mex}
\big|U(t)\big| _{Y^s} \leq \big|U_0\big|_{Y^s}   e^{ \tilde{\lambda}\max(\eps,\beta)t}  + \max(\eps,\beta) \int^{t}_{0}  \tilde{C} \big( h_{min}^{-1},q^{-1}_{min}, \vert b \vert _{H^{s+3}}  , |\underline{U} (t')|_{Y^s} \big) \;  e^{ \tilde{\lambda} \max(\eps,\beta) ( t-t')}  \;dt'  \; ,
\end{equation}
\begin{equation}\label{derivative-energy-est-Mex}
\vert \partial_tU (t)\vert _{Y^{s-1}}   \le  \tilde{\tilde{C}} |U(t) |_{Y^s} + \beta C\big(h_{min}^{-1}, \vert b  \vert_{H^{s}}, \vert \underline{U}\vert_{Y^s}\big) \;,
\end{equation}
for some positive constant $\tilde{\lambda} >0$ depending on $M_{CH}$, $h_{min}^{-1}$, $q^{-1}_{min}$, and the norms $ \vert b  \vert _{H^{s+2}}$, $\vert\partial_t \underline{\zeta}\vert_{L^{\infty}}$, $\vert\partial_t \underline{v}\vert_{L^{\infty}}$ and $\vert \underline{U}\vert_{Y^s}$, where $\max(\eps,\beta) \tilde{C}$ is the upper bound of the inner product $(\Lambda^s B(\underline{U}),Z[\underline U]\Lambda^s v)$ and $\tilde{\tilde{C}}$ a positive constant depending also on $h_{min}^{-1}$, $q^{-1}_{min}$ and the norms $ \vert b \vert _{H^{s+3}}$ and $\vert \underline{U}\vert_{Y^s}$.

\end{proposition}

\begin{proof}

Once the $Y^s$ energy estimate \eqref{energy-est-Mex} is obtained, the existence, uniqueness, and regularity of the solution to the linearized system \eqref{SSlsys} is determined using the Cauchy-Lipschitz technique, which also necessitates a similar estimate. As a result, we do not provide details on the implementation strategy because the majority of the work is to derive a prior energy estimate (for similar details, see \cite{DucheneIsrawiTalhouk14,IsrawiLteifTalhouk15,LteifIsrawi18}). However, we will concentrate on demonstrating the key steps \eqref{energy-est-Mex}-\eqref{derivative-energy-est-Mex}.

Let us consider any $\lambda\in\RR$ to be fixed later. Since the symmetrizer $Z[\underline{U}]$ is symmetric, it is not hard to check that we have
\begin{multline}\label{estimations}
\frac{1}{2}e^{\max(\eps,\beta )  \lambda t}\partial_t(e^{- \max(\eps,\beta )  \lambda t}E^s(U)^2)=
- \lambda  \frac{\max(\eps,\beta ) }{2}E^s(U)^2 
-\big(A[\underline{U}]\Lambda^s\partial_x U, Z[\underline U]\Lambda^s U\big) \\- \big(\big[\Lambda^s,A[\underline{U}]\big] \partial_xU,Z[\underline U]\Lambda^s U\big)
-\big(\Lambda^sB(\underline{U}),Z[\underline U]\Lambda^sU\big)+\frac{1}{2}\big(\Lambda^s\zeta, \big[\partial_t,  \dfrac{\underline{Q}_0+\epsilon^2 \underline{Q}_1}{H(\epsilon\underline{\zeta},\beta b)} \big ]\Lambda^s\zeta\big)+\frac{1}{2}\big(\Lambda^sv,[\partial_t,\underline{\mfT}]\Lambda^sv\big) \; .
\end{multline}
The key point is to bound from above each term at the right-hand side of \eqref{estimations} in terms of $E^s(U)$ and $E^s(\underline{U})$, then for a specific choice of $\lambda$ and by Gr$\ddot{\text{o}}$nwall's inequality. The prior energy estimates \eqref{energy-est-Mex}-\eqref{derivative-energy-est-Mex} follows.\\

$\bullet$ \underline{Estimate of $\big(Z[\underline U]A[\underline U]\partial_x \Lambda^s U, \Lambda^s U\big)$.} By definition, it holds that
 \begin{equation*}
Z[\underline{U}]A[\underline{U}]=\begin{pmatrix}
 \epsilon \dfrac{\underline{Q}_0+\epsilon^2 \underline{Q}_1}{H(\epsilon\underline{\zeta},\beta b)} H'(\epsilon\underline{\zeta},\beta b) \underline{v}& \underline{Q}_0+\epsilon^2 \underline{Q}_1 \\
\underline{Q}_0+\epsilon^2 \underline{Q}_1   & \mfQ[\epsilon\underline{\zeta},\beta b,\underline{v}]+\epsilon \mfT[\epsilon\underline{\zeta},\beta b]( \varsigma\underline{v}\cdot)
\end{pmatrix}.
\end{equation*}
Therefore, we have
\begin{align*}
\big(Z[\underline{U}]A[\underline{U}]\partial_x \Lambda^s  U, \Lambda^s  U\big)&=\Big(\epsilon\dfrac{\underline{Q}_0+\epsilon^2 \underline{Q}_1}{H(\epsilon\zetabar,\beta b)}H'(\epsilon\zetabar,\beta b)\vbar\partial_x \Lambda^s \zeta, \Lambda^s\zeta \Big)+ \Big(Q_0(\epsilon\zetabar,\beta b)\partial_x \Lambda^s v,\Lambda^s\zeta\Big)+\Big(\epsilon^2Q_1(\epsilon\zetabar,\beta b,\vbar)\partial_x\Lambda^s v,\Lambda^s\zeta \Big)\\ 
& \quad+\Big(Q_0(\epsilon\zetabar,\beta b)\partial_x\Lambda^s \zeta, \Lambda^sv \Big) +\Big(\epsilon^2Q_1(\epsilon\zetabar,\beta b,\vbar)\partial_x \Lambda^s\zeta, \Lambda^sv \Big) +\Big( \mfQ[\epsilon\underline{\zeta},\beta b,\underline{v}] \partial_x \Lambda^sv, \Lambda^sv \Big) \\
&\quad + \epsilon \Big( \mfT[\epsilon\underline{\zeta},\beta b](\varsigma\underline{v}\partial_x\Lambda^s v), \Lambda^sv \Big) = I_1+..+I_7 \; .
\end{align*}
To control $I_1$, by integration by parts with Lemma \ref{lemma4} in hands and the continuous embedding $H^s(\RR)\subset W^{1,\infty}(\RR)$ for $s>3/2$, it holds that
\begin{align*}
I_1 &= -\dfrac{1}{2}\Big(\epsilon\partial_x\big(\dfrac{Q_0(\epsilon\zetabar,\beta b)+\epsilon^2Q_1(\epsilon\zetabar,\beta b,\vbar)}{H(\epsilon\zetabar,\beta b)}H'(\epsilon\zetabar,\beta b)\vbar\big) \Lambda^s\zeta, \Lambda^s\zeta \Big)\\
&  \le \eps C\big( h_{min}^{-1},q_{min}^{-1}, \vert b\vert _{H^{s+1}}, \vert \underline{U}\vert_{Y^s}\big) E^s(U)^2 \; .
\end{align*}
To control $(I_2+I_4)+(I_3+I_5)$, by integration by parts with Lemma \ref{lemma4} in hands, it holds that
\begin{align*}
(I_2+I_4)+(I_3+I_5) & = -\Big(\partial_x\big(Q_0(\epsilon\zetabar,\beta b)\big)  \Lambda^s\zeta , \Lambda^sv \Big)-\epsilon^2\Big(\partial_x\big(Q_1(\epsilon\zetabar,\beta b,\vbar)\big) \Lambda^s\zeta,  \Lambda^sv\Big) \\
&\le \max(\eps,\beta) C\big( h_{min}^{-1},q_{min}^{-1}, \vert b\vert _{H^{s+2}}, \vert \underline{U}\vert_{Y^s}\big) E^s(U)^2 \; .
\end{align*}
To control $I_6$, we recall the expression \eqref{defmfQ}. By integration by parts with Lemma \ref{lemma4} in hands, it holds that
\begin{align*}
I_6=&-\frac{1}{2} \Big(\ \epsilon\partial_x\big[ q_1(\epsilon\underline{\zeta},\beta b)(H'(\epsilon\zetabar,\beta b)-\varsigma)\vbar\big] \Lambda^sv, \Lambda^sv\Big)-\mu\Big([\mathcal{A}] \vbar \Lambda^sv_x , \Lambda^sv\Big)  -\mu\Big([\mathcal{B}]  \vbar_x\Lambda^sv_x , \Lambda^sv\Big)\\
& +\mu\Big(\partial_x\big([\mathcal{C}] \vbar\big) \Lambda^s v_x  , \Lambda^sv\Big)+\mu\Big( [\mathcal{C}] \vbar \Lambda^s v_x  , \Lambda^sv_x\Big) +\mu\Big(\partial_x\big([\mathcal{D}]\big) \vbar_x \Lambda^sv_x,\Lambda^sv \Big) +\mu\Big( [\mathcal{D}]  \vbar_x \Lambda^sv_x,\Lambda^sv_v \Big)\\
&\le \max(\eps,\beta) C\big( h_{min}^{-1},q_{min}^{-1}, \vert b\vert _{H^{s+3}}, \vert \underline{U}\vert_{Y^s}\big) E^s(U)^2 \; .
\end{align*}
To control $I_7$, we recall the expression \eqref{defZ}. By integration by parts with Lemma \ref{lemma4} in hands, it holds that
\begin{align*}
I_7 & = -\frac12\Big(\ \partial_x(q_1(\epsilon\underline{\zeta},\beta b)\varsigma\underline{v}) \Lambda^sv,\Lambda^sv \Big)+\mu\epsilon\beta\Big((\kappa_0\partial_x \zetabar\partial_x b \varsigma\vbar) \Lambda^sv_x , \Lambda^sv \Big) +\mu\Big(\ \nu q_2(\epsilon\underline{\zeta},\beta b)\partial_x(\varsigma\underline{v} \Lambda^sv_x ) , \Lambda^sv_x\ \Big)\\ 
& =-\frac12\Big(\ \partial_x(q_1(\epsilon\underline{\zeta},\beta b)\varsigma\underline{v}) \Lambda^sv,\Lambda^sv \Big)+\mu\epsilon\beta\Big((\kappa_0\partial_x \zetabar\partial_x b \varsigma\vbar) \Lambda^sv_x , \Lambda^sv \Big)   +\mu\Big( \nu q_2(\epsilon\underline{\zeta},\beta b)(\partial_x\varsigma \underline{v} ) \Lambda^sv_x  ,  \Lambda^sv_x \Big)  \\
& \quad -\mu\frac12\Big( \partial_x(\nu q_2(\epsilon\underline{\zeta},\beta b)\varsigma\underline{v} ) \Lambda^sv_x  ,  \Lambda^sv_x \Big)\\
&\le \max(\eps,\beta) C\big( \vert b\vert _{H^{s+1}},q_{min}^{-1}, \vert \underline{U}\vert_{Y^s}\big) E^s(U)^2 \; .
\end{align*}
As a conclusion from the above estimates, it holds that
\begin{equation}\label{est1}
| \big(Z[\underline U]A[\underline U]\partial_x \Lambda^s U, \Lambda^s U\big) | \le  \max(\eps,\beta) C\big( h_{min}^{-1}, \vert b\vert _{H^{s+3}}, \vert \underline{U}\vert_{Y^s}\big) E^s(U)^2 \; .
\end{equation}

$\bullet$ \underline{Estimate of $ \big(\big[\Lambda^s,A[\underline{U}]\big]\partial_x  U,Z[\underline{U}]\Lambda^s U\big)$}. Using the definition of $A[\cdot]$ and $Z[\cdot]$ in~\eqref{defA0A1} and \eqref{defZ}, one has
 \begin{align*}
\big(\big[\Lambda^s,A[\underline{U}]\big]\partial_x U,Z[\underline{U}]\Lambda^s U\big)&= \epsilon\Big([\Lambda^s, H'(\epsilon\underline{\zeta},\beta b) \underline{v}] \zeta_x, \dfrac{\underline{Q}_0+\epsilon^2 \underline{Q}_1}{H(\epsilon\underline{\zeta},\beta b)} \Lambda^s\zeta\Big)  + \Big( [\Lambda^s, H(\epsilon\underline{\zeta},\beta b)] v_x, \dfrac{\underline{Q}_0+\epsilon^2 \underline{Q}_1}{H(\epsilon\underline{\zeta},\beta b)} \Lambda^s\zeta\Big) \\
& +\Big([\Lambda^s, \underline{\mfT}^{-1}\big(\underline{Q}_0\cdot+\epsilon^2 \underline{Q}_1\cdot\big)]  \zeta_x , \underline{\mfT}\Lambda^s v\Big) +\Big([\Lambda^s, \underline{\mfT}^{-1}(\mfQ[\epsilon\underline{\zeta},\beta b,\underline{v}] \cdot)  ] v_x, \underline{\mfT}\Lambda^s v\Big) \\
& + \Big([\Lambda^s , \epsilon\varsigma \underline v] v_x, \underline{\mfT}\Lambda^s v\Big)  \equiv  J_1 + .. + J_5.
\end{align*}
To control of $J_1$, using \eqref{cest} with Lemma \ref{lemma4} in hands and the fact that $| \partial_x(H'(\eps\underline{\zeta},\beta b )\underline{v}) |_{H^{s-1}}  \le C(h_{min}^{-1},|\underline{\zeta}|_{H^{s}},,\vert b\vert _{H^s})|\underline{v}|_{H^{s}}$, it holds that
$$
J_1 \le \eps C\big(  h_{min}^{-1},q_{min}^{-1}, \vert b\vert _{H^{s}}, \vert \underline{U}\vert_{Y^s}\big) E^s(U)^2 \;.
$$
Similarly, to control $J_2$, since $|\partial_x H(\eps\underline{\zeta},\beta b) |_{H^{s-1}} \le   \max(\eps,\beta)C(h_{min}^{-1},|\underline{\zeta}|_{H^{s} },\vert b\vert _{H^s})$, it holds that
$$
J_2 \le \max( \eps,\beta) C\big(  h_{min}^{-1},q_{min}^{-1}, \vert b\vert _{H^{s}}, \vert \underline{U}\vert_{Y^s}\big) E^s(U)^2 \;.
$$
Now we introduce the following useful commutator identities for the rest of the terms
\begin{equation}\label{commu}
[\Lambda^s,\partial_x(f\partial_x\cdot)]g = \partial_x[\Lambda^s,f]g_x, \qquad \qquad 
[\Lambda^s,f\partial_x]g=\partial_x[\Lambda^s,f]g-[\Lambda^s,\partial_xf]g=[\Lambda^s,f]\partial_xg \; .
\end{equation}
To control $J_3$, remark that $\mfT$ is symmetric with
\begin{equation*}
\underline{\mfT}[\Lambda^s, \underline{\mfT}^{-1}] (\underline{Q}_0\cdot+\epsilon^2 \underline{Q}_1\cdot )\zeta_x=\underline{\mfT}[\Lambda^s, \underline{\mfT}^{-1}(\underline{Q}_0\cdot+\epsilon^2 \underline{Q}_1\cdot)]\zeta_x-[\Lambda^s, \underline{Q}_0\cdot+\epsilon^2 \underline{Q}_1\cdot ]\zeta_x \; .
\end{equation*}
Moreover, since $[\Lambda^s,\underline{\mfT}^{-1}]=-\underline{\mfT}^{-1}[\Lambda^s,\underline{\mfT}]\underline{\mfT}^{-1}$, one gets
\begin{equation*}
\underline{\mfT}[\Lambda^s, \underline{\mfT}^{-1}\;(\underline{Q}_0\cdot+\epsilon^2 \underline{Q}_1\cdot)]\zeta_x=-[\Lambda^s, \underline{\mfT}] \underline{\mfT}^{-1}\big(\underline{Q}_0\zeta_x+\epsilon^2 \underline{Q}_1\zeta_x\big) +[\Lambda^s,\underline{Q}_0 +\epsilon^2 \underline{Q}_1]\zeta_x \; .
\end{equation*}
Therefore, we can now use Corollary~\ref{col:comwithT} in view of expression \eqref{defQ0Q1}, \eqref{commu}, such that we have
\begin{align*}
J_3 &=  -\big(\big[\Lambda^s, \underline{\mfT}\big] \underline{\mfT}^{-1}(\underline{Q}_0\cdot+\epsilon^2 \underline{Q}_1\cdot) \zeta_x , \Lambda^s v\big) + \big(  [\Lambda^s,\underline{Q}_0 +\epsilon^2 \underline{Q}_1]  \zeta_x, \Lambda^s v\big) \\
& \le  \max( \eps,\beta) C\big(M_{CH},  h_{min}^{-1},q_{min}^{-1}, \vert b\vert _{H^{s+2}}, \vert \underline{U}\vert_{Y^s}\big) E^s(U)^2 \;.
\end{align*}
Similarly, recall that $\mfQ[\epsilon\underline{\zeta},\beta b,\underline{v}]f =\Big(\epsilon q_1(\epsilon\underline{\zeta},\beta b)(H'(\epsilon\underline{\zeta},\beta b) -\varsigma)-\mu[\mathcal{A}]\Big)\underline{v}f-\mu[\mathcal{B}]\underline{v}_x f -\mu[\mathcal{C}]\underline{v}\partial_x f-\mu[\mathcal{D}]\partial_x(\underline{v}_x f)$. In view of Lemma \ref{lemma4} and \eqref{commu}, it holds that
\begin{align*}
J_4 &=  -\big( \big[\Lambda^s, \underline{\mfT}\big] \underline{\mfT}^{-1}(\mfQ[\epsilon\underline{\zeta},\beta b,\underline{v}]\cdot) \zeta_x , \Lambda^s v\big) + \big(  [\Lambda^s,\mfQ[\epsilon\underline{\zeta},\beta b,\underline{v}]]  \zeta_x, \Lambda^s v\big) \\
& \le  \max( \eps,\beta) C\big( M_{CH},h_{min}^{-1},q_{min}^{-1}, \vert b\vert _{H^{s+3}}, \vert \underline{U}\vert_{Y^s}\big) E^s(U)^2 \;.
\end{align*}
To control $J_5$, using \eqref{commu} and the expression of \eqref{defZ} and by integration by parts, it holds that
\begin{align*}
J_5 &=  \big([\Lambda^s , \epsilon\varsigma \underline v] v_x,  q_1(\eps\underline{\zeta},\beta b ) \Lambda^s v\big) +  \eps\mu \big([\Lambda^s , \epsilon\varsigma \underline v] v_x,  \kappa_0\underline{\zeta}_x b_x  \Lambda^s v\big)
+ \mu \big([\Lambda^s , \epsilon\varsigma \underline v] v_{xx},  \nu q_2(\eps\underline{\zeta},\beta b ) \Lambda^s v_x\big)\\
&\quad +  \mu\big([\Lambda^s , \epsilon\partial_x(\varsigma \underline v)] v_x,  \nu q_2(\eps\underline{\zeta},\beta b ) \Lambda^s v_x\big)\\
&  \le  \max( \eps,\beta) C\big( h_{min}^{-1},q_{min}^{-1}, \vert b\vert _{H^{s+1}}, \vert \underline{U}\vert_{Y^s}\big) E^s(U)^2 \;.
\end{align*}
As a conclusion from the above estimates, it holds that
\begin{equation}\label{est:2}
|\big(\big[\Lambda^s,A[\underline{U}]\big]\partial_x U,S[\underline{U}]\big]\Lambda^s U\big) |  \max( \eps,\beta) C\big( M_{CH},h_{min}^{-1},q_{min}^{-1}, \vert b\vert _{H^{s+3}}, \vert \underline{U}\vert_{Y^s}\big) E^s(U)^2 \;.
\end{equation}

$\bullet$ \underline{Estimate of  $\big(\Lambda^s B[\underline U],Z[\underline U] \Lambda^s  U\big)$.} Recalling the expressions \eqref{defB} and \eqref{defZ}, it is not hard to check that
\begin{align*}
\big(\Lambda^s B[\underline U],Z[\underline U] \Lambda^s  U\big) & =\frac{1}{2} \Big( \Lambda^s\big(- \beta b_x G(\epsilon\underline{\zeta},\beta b)\underline{v} \big),  \dfrac{\underline{Q}_0+\epsilon^2 \underline{Q}_1}{H(\epsilon\underline{\zeta},\beta b)} \Lambda^s\zeta\Big) \\
&\quad +\Big( \big [\Lambda^s,\underline{\mfT}^{-1} \big]\big(\dfrac{\gamma\epsilon\beta q_1(\epsilon\zetabar,\beta b)\underline{h}_1(\underline{h}_1+\underline{h}_2)\underline{v}^2\partial_x b}{(\underline{h}_1+\gamma \underline{h}_2)^3}-\mu[\mathcal{E}] \vbar^2 \big), \underline{\mfT}\Lambda^s v\Big) \\
&\quad +\Big( \Lambda^s\big(\dfrac{\gamma\epsilon\beta q_1(\epsilon\zetabar,\beta b)\underline{h}_1(\underline{h}_1+\underline{h}_2)\underline{v}^2\partial_x b}{(\underline{h}_1+\gamma \underline{h}_2)^3}-\mu[\mathcal{E}] \vbar^2 \big),  \Lambda^s v\Big) = B_1+B_2+B_3.
\end{align*}
 It is not hard to check that $B_1$ and $B_3$ can be controlled by $\max(\eps,\beta ) C\big( h_{min}^{-1},q_{min}^{-1}, \vert b\vert _{H^{s+3}}, \vert \underline{U}\vert_{Y^s}\big) E^s(U)$. However, as for $J_3$, one may control $B_2$ by $\max(\eps,\beta ) C\big(M_{CH}, h_{min}^{-1},q_{min}^{-1}, \vert b\vert _{H^{s+3}}, \vert \underline{U}\vert_{Y^s}\big) E^s(U)$. Hence, it holds that
 \begin{equation}\label{est2}
 \big(\Lambda^s B[\underline U],Z[\underline U] \Lambda^s  U\big) \le \max(\eps,\beta ) C\big(M_{CH}, h_{min}^{-1},q_{min}^{-1}, \vert b\vert _{H^{s+3}}, \vert \underline{U}\vert_{Y^s}\big) E^s(U) \; .
 \end{equation}

$\bullet$ \underline{Estimate of $\frac{1}{2}\big(\Lambda^s\zeta, \big[\partial_t,  \dfrac{\underline{Q}_0+\epsilon^2 \underline{Q}_1}{H(\epsilon\underline{\zeta},\beta b)} \big ]\Lambda^s\zeta\big)$.}
Since $[\partial_t, f] g = f_tg$ and the independent of the bottom profile on the spatial dimension, it holds that
\begin{align}\label{est3}
\frac{1}{2}\big(\Lambda^s\zeta, \big[\partial_t,  \dfrac{\underline{Q}_0+\epsilon^2 \underline{Q}_1}{H(\epsilon\underline{\zeta},\beta b)} \big ]\Lambda^s\zeta\big) &= \Big(\Lambda^s\zeta,\partial_t \Big(\frac{Q_0(\epsilon\underline{\zeta},\beta b)+\epsilon^2 Q_1(\epsilon\underline{\zeta},\beta b,\underline{v})}{H(\epsilon\underline{\zeta},\beta b)}\Big)\Lambda^s\zeta\Big) \nonumber \\
& \le \max(\eps,\beta ) C\big( h_{min}^{-1},q_{min}^{-1}, \vert b\vert _{H^{s+2}}, \vert\partial_t \underline{\zeta}\vert_{L^{\infty}},\vert\partial_t \underline{v}\vert_{L^{\infty}}\big)E^s(U)^2 \; .
\end{align}

$\bullet$ \underline{Estimate of $\frac{1}{2}(\Lambda^sv,\big[\partial_t,\underline{\mfT}\big]\Lambda^sv)$.}
As above, using the expression of \eqref{defZ}, by integration by parts, it holds that
\begin{align}\label{est4}
 (\Lambda^sv,\big[\partial_t,\underline{\mfT}\big]\Lambda^sv) & =  \epsilon\Big(\Lambda^s v,\kappa_1(\partial_t \underline{\zeta}) \Lambda^s v\Big)- 2\mu\epsilon\beta\Big(\Lambda^s v,\kappa_0 \partial_t \zetabar \partial_x b \ \Lambda^s   v_x\Big)  -\mu\epsilon\beta\Big(\Lambda^s v,\partial_x \kappa_0 \partial_t \zetabar \partial_x b \ \Lambda^s v\Big) \nonumber \\
 & \quad -\mu\epsilon\beta\Big(\Lambda^s v,\kappa_0 \partial_t \zetabar \partial_x^2 b \ \Lambda^s v\Big) +\mu\epsilon\Big(\Lambda^s \partial_x v,\nu\kappa_2(\partial_t\underline{\zeta})\Lambda^s \partial_x v \Big) \nonumber \\
 &\le \max(\eps,\beta ) C\big( h_{min}^{-1},q_{min}^{-1}, \vert b\vert _{H^{s+2}}, \vert\partial_t \underline{\zeta}\vert_{L^{\infty}}\big) E^s(U)^2 \; .
\end{align}
Thanks to the above estimates \eqref{est1}-\eqref{est:2}-\eqref{est2}-\eqref{est3}-\eqref{est4}, combine these with \eqref{estimations} it holds that 
\begin{align*}
\frac{1}{2}e^{  \lambda   \max(\eps,\beta )   t} \partial_t (e^{- \lambda   \max(\eps,\beta )   t}E^s(U)^2) & \le   \max(\eps,\beta ) \big( C\big( M_{CH},h_{min}^{-1},q_{min}^{-1}, \vert b\vert _{H^{s+2}}, \vert\partial_t \underline{\zeta}\vert_{L^{\infty}},\vert\partial_t \underline{v}\vert_{L^{\infty}}\big)  - \frac{\lambda}{2} \big)E^s(U)^2\\
 &\quad +  \max(\eps,\beta ) C\big( h_{min}^{-1},q_{min}^{-1}, \vert b\vert _{H^{s+3}}, \vert \underline{U}\vert_{Y^s}\big) E^s(U)\; .
\end{align*}
Now, for all $0\le \max(\eps,\beta ) t\le T$, we take $\tilde{\lambda} = \lambda\ge 2 C\big( M_{CH}, h_{min}^{-1},q_{min}^{-1}, \vert b\vert _{H^{s+2}}, \vert\partial_t \underline{\zeta}\vert_{L^{\infty}},\vert\partial_t \underline{v}\vert_{L^{\infty}}\big)>0$ so that the differential inequality below holds:
$$
\frac{d}{dt}  E^s(U)  \le \frac{1}{2}\tilde{\lambda}\max(\eps,\beta )E^s(U) +  \max(\eps,\beta ) C\big( h_{min}^{-1},q_{min}^{-1}, \vert b\vert _{H^{s+3}}, \vert \underline{U}\vert_{Y^s}\big)  \; .
$$
At this stage, we multiply the above differential inequality by $1/\sqrt{ e^{ \tilde{\lambda} \max(\eps,\beta )  ( t-t')} } $ then we integrate on $(0,t)$  for all $0\le \max(\eps,\beta ) t\le T$, so that the desired energy estimate \eqref{energy-est-Mex} holds.

Now to establish the derivative energy estimate \eqref{derivative-energy-est-Mex} we use the linearized system \eqref{SSlsys} with \eqref{defA0A1}-\eqref{defB}. Indeed, by definition we have
$$
\big\vert\partial_t U\big\vert_{Y^{s-1}} ^2 = \big\vert -A[\underline{U}]\partial_x U-B (\underline{U})\big\vert_{Y^{s-1}}^2   =  |\mathfrak{A}|_{H^{s-1}}^2+ |\mathfrak{B}|_{H^{s-1}}^2+ \mu|\partial_x\mathfrak{B}|_{H^{s-1}}^2
$$
where $\mathfrak{A} = \epsilon H'(\epsilon \underline{ \zeta},\beta b) \underline{v}\zeta_x + H(\epsilon\underline{\zeta},\beta b)v_x  -\beta b_xG(\epsilon\underline{\zeta},\beta b)\underline{v}$ and
$$
\mathfrak{B} =\epsilon\varsigma \underline{v}v_x+ \underline{\mfT}^{-1}\Big(\underline{Q}_0 \zeta_x+\epsilon^2\underline{Q}_1\zeta_x+  \mfQ[\epsilon\underline{\zeta},\beta b,\underline{v}] \zeta_x  + \gamma\epsilon\beta q_1(\epsilon\underline{\zeta},\beta b)\underline{h}_1(\underline{h}_1+\underline{h}_2)(\underline{h}_1+\gamma \underline{h}_2)^{-3} b_x \underline{v}^2-\mu[\mathcal{E}]\underline{v}^2\Big).
 $$
With Lemma \ref{lemma4} in hands and using \eqref{i-1}, it holds that
\begin{equation}\label{utest}
 \big\vert\partial_t U\big\vert_{Y^{s-1}} \le C\big(h_{min}^{-1}, \vert b  \vert_{H^{s+3}}, \vert \underline{U}\vert_{Y^s}\big)\vert U\vert_{Y^s} + \beta C\big(h_{min}^{-1}, \vert b  \vert_{H^{s}}, \vert \underline{U}\vert_{Y^s}\big)\;.
\end{equation}

It remains to show that $U$ satisfies preserves the depth condition. Consequently, using \eqref{utest} and the assumption \eqref{CondDepth}, it holds that for all $0\le \max(\eps,\beta) t\le T$ we have
\begin{align*}
h_{1,2} & =h_{1,2}(0,x) + \eps\int_{0}^t \partial_{t'} \zeta(t',x)\; dt'   \\
& \ge  h_{min} - \eps    \Big( C_2(h_{min}^{-1}, \vert b\vert_{H^{s+3}},\vert \underline{U}\vert_{Y^s},\vert \underline{U}_t\vert_{Y^s},\vert U_0\vert_{Y^s})  + C_3(h_{min}^{-1}, \vert b\vert_{H^{s+3}},\vert \underline{U}\vert_{Y^s})\; e^{T} \Big) T \\
& \ge  h_{min} -     \Big( C_2(h_{min}^{-1}, \vert b\vert_{H^{s+3}},\vert \underline{U}\vert_{Y^s},\vert \underline{U}_t\vert_{Y^s},\vert U_0\vert_{Y^s})  T  + C_3(h_{min}^{-1}, \vert b\vert_{H^{s+3}},\vert \underline{U}\vert_{Y^s} )  \; e^{2T} \Big)  \;.
\end{align*}
Therefore it is possible to choose small enough time less than or equal
$
\min \left\{ \frac{h_{min}}{4 C_2} , \frac{1}{2}\ln \left( \frac{h_{min}}{4 C_3} \right) \right\} \; .
$
 Thus $U$ satisfies the depth-condition \eqref{CondDepth} with $h_{min}$ replaced by its one half.

\end{proof}

\subsection{Well-posedness}
\begin{theorem}[Local existence]\label{theom1}
Let $s>3/2$. Also let $b\in H^{s+3}(\RR)$ and an initial data $U_0= (\zeta_0,v_0)^T\in Y^{s}$ endowed with the norm \eqref{endowed-norm}
such that for any $x\in\RR$ 
 the depth condition~\eqref{CondDepth} is satisfied in addition to the assumptions~\eqref{condbot},~\eqref{CondEllipticity} and~\eqref{H3}.
 Then there exists a lifetime $T_{mGN}>0$ and a unique solution $U\in C\big([0,\frac{T_{mGN}}{\max(\eps,\beta)}],Y^s\big)$ to model \eqref{condensedeq}, which depends continuously on the initial conditions $\vert U_0\vert_{Y^s}$ and $h_{min}^{-1},q_{min}^{-1}>0$. Additionally, we have the solution size estimate
$$
\vert U\vert _{Y^s} + \vert \partial_tU\vert _{Y^{s-1}} \le  C\big( M_{CH},h_{min}^{-1},q_{min}^{-1}, \vert b\vert _{H^{s+3}}, \vert U_0\vert_{Y^s}\big) ,\qquad for \quad   0\le \max(\eps,\beta )t\le  T_{mGN} \; .
$$
Also, the depth condition \eqref{CondDepth} is still satisfied on the existence-interval $[0,\frac{T_{mGN}}{\max(\eps,\beta)}]$.
\end{theorem}
\begin{proof}
With Proposition \ref{prop1}, i.e. the mathematical analysis of the linearized system, the local well-posedness result on the interval $[0,\max(\eps,\beta)^{-1}T_{mGN}]$ is obtained by an adaptation of the standard well-posedness proof of hyperbolic systems (see Theorem 7.3 in \cite{DucheneIsrawiTalhouk14} or Theorem 6.1 in \cite{IsrawiLteifTalhouk15} ).\\ 
We will focus on the key points of the proof. The proof techniques are based on those for hyperbolic systems inspired from Chapter III B.1 of \cite{AlinhacGerard}. We will give a brief description of the proof, which is initiated by a series of nonlinear problems via the induction relation
\begin{equation}\label{UN}
\forall\hspace{0.1cm}n \in\mathbb{N}\hspace{0.1cm},\qquad
\left\{
\begin{array}{lcl}
\displaystyle \partial_tU^{n+1} + A[U^{n}] U^{n+1} + B(U^{n}) = 0 , \\
\displaystyle U^{n+1}_{\mid_{t = 0}} = U_0, \qquad \text{ with }\qquad U^0 = U_0 \; .
\end{array}
\right.
\end{equation}
From the estimates \eqref{energy-est-Mex}-\eqref{derivative-energy-est-Mex} of section \ref{linear-analysis}, the convergence of solution $U^n=(\zeta^n,v^n)$ for system \eqref{UN}, can be established combined with standard arguments. Indeed, by induction on $n$ it holds that $U^{n+1} \in C\big([0,\frac{T_{mGN}}{\max(\eps,\beta)}] ,Y^s\big)$ such that $\vert U^{n+1}\vert_{Y^s}\lesssim \vert U_0\vert_{Y^s}$ and $\vert \partial_tU^{n+1}\vert_{Y^{s-1}}\lesssim \vert U_0\vert_{Y^s}$ for all $t\in[0,\frac{T_{mGN}}{\max(\eps,\beta)}]$. 
 Consequently, $\tilde{\lambda}$  does not depend on $\vert \partial_tU\vert_{L^{\infty}}$ anymore. As a result, the depth-condition will be satisfied by $U^{n+1}(t)$ for small enough $T(\vert U_0\vert_{Y^s})$. It is still necessary to demonstrate the convergence of $U^n$ towards \eqref{condensedeq} solution. In fact, the difference between two consecutive approximate solutions $V^n=U^{n+1}-U^n$ reads the system
 \begin{equation}\label{VN}
\forall\hspace{0.1cm}n \in\mathbb{N}\hspace{0.1cm},\qquad
\left\{
\begin{array}{lcl}
\displaystyle \partial_tV^{n} + A[U^n] V^{n} = -\big(A[U^n]-A[U^{n-1}]\big) U^n + B(U^{n}) - B(U^{n-1})  ,\\
\displaystyle V^{n}_{\mid_{t = 0}} = 0 .
\end{array}
\right.
\end{equation}
In the case $s=0$, one may check that the energy estimate \eqref{energy-est-Mex} for $V^n$ on $[0,\max(\eps,\beta)^{-1} T]$ become
$$
\vert V^n(t)\vert_{Y^0} \le \frac{1}{n!} \max(\eps,\beta)^n t^nC^n(\vert U^n\vert_{W^{1,\infty}}) \sup _{t'\in[0,\max(\eps,\beta)^{-1}T_{mGN}]}\vert V^0(t')\vert _{Y^0} \; ,
$$
and thus the sequence $U^n=U_0+\sum_0^{n-1}V^i$ converges in $Y^0_{T_{mGN}}$. Consequently, an interpolation argument suggests that $U^n\in Y^{s'}_{T_{mGN}}$ for $s'<s$. As a result, selecting $s'-1>1/2$, the limit $U\in Y^s_{T_{mGN}}$ of the iterative scheme \eqref{UN} is a unique solution of system \eqref{condensedeq} and satisfies the energy estimate \eqref{energy-est-Mex}. Finally, the analysis above shows that the maximal lifetime $T_{\max}$ is bounded from below by some $\max(\eps,\beta)^{-1}T_{mGN}>0$, whereas if $T_{\max}<\infty$, the behavior of the solution as $t \rightarrow T_{\max}$ follows from standard continuation arguments.

\end{proof}

\begin{remark}
Our model's full justification (consistency+convergence) follows in the same way as the models studied in   \cite{DucheneIsrawiTalhouk14,IsrawiLteifTalhouk15,LteifIsrawi18}.
\end{remark}

  \section{Surface water-wave equations (One-layer case)}\label{Sec7}
We conclude our work in this section by using the previous section's results to derive some new results on the propagation of surface water-waves (one-layer case with free surface) over uneven bottoms while accounting for surface tension effects. In fact, in the water-wave case, the GNCH with uneven bottoms can be recovered by simply setting $\gamma=0$ and $\delta=1$ in~\eqref{eq:Serre2}. Thus, one gets
      \begin{equation}\label{eq:Serre1}\left\{ \begin{array}{l}
      \partial_{ t}\zeta +\partial_x\big( h v\big)\ =\  0,\\ \\
      \mathfrak T[\epsilon\zeta,\beta b] \left( \partial_{ t}   v +  \epsilon \varsigma { v } \partial_x {  v } \right) +q_1(\epsilon\zeta,\beta b)\partial_x
      \zeta \\ =    \mu\Big([\mathcal{A}]v\partial_x v +  [\mathcal{B}](\partial_x v)^2 +[\mathcal{D}]\partial_x\big((\partial_x v)^2\big)+[\mathcal{E}]v^2+[\mathcal{F}]\partial_x \zeta \Big),
      \end{array} \right.  \end{equation}  
where we denote $h= 1+\epsilon\zeta-\beta b$ and 
 \begin{equation*}\label{defT2}
      {\mathfrak T}[\epsilon\zeta,\beta b]V \ = \ q_1(\epsilon\zeta,\beta b)V \ + \ \mu\epsilon\beta\kappa_0 \partial_x \zeta \partial_x b V \ - \ \mu\nu \partial_x\Big(q_2(\epsilon\zeta,\beta b)\partial_xV \Big),
      \end{equation*}
      with $q_1(X,Y)\equiv 1+\kappa_1 X +(\omega_1 + \epsilon \phi_1) Y $ and $q_2(X,Y)\equiv 1+\kappa_2 X +\omega_2 Y$. One can easily check that in the one-layer case ($\gamma=0$ and $\delta=1$), the functions $\nu, \kappa_0, \kappa_1, \kappa_2, \phi_1, \omega_1, \omega_2$ and $\varsigma$ defined respectively in~\eqref{defnu},~\eqref{defkappa0},~\eqref{defkappa1},~\eqref{defkappa2},~\eqref{defphi1},~\eqref{ode2},~\eqref{ode1} and~\eqref{defvarsigma} become:
         \begin{equation}\label{defnu2}
 \nu \ = \ \frac13-\frac1{\bo},
      \end{equation}
 \begin{equation}\label{defkappa02}
 \kappa_0 \ = \ 1+\beta\omega_1 b, 
  \end{equation}
\begin{equation}\label{defkappa12}
 \kappa_1 \ = \ \dfrac{(1+\beta \omega_1 b)(1-\beta b)}{(1-\beta b)^2-\dfrac{3}{\bo}}, 
  \end{equation}
  \begin{equation}\label{defkappa22}
\kappa_2 \ = \ \dfrac{(1-\beta b)(1+\beta \omega_1 b)}{\nu}, \end{equation}
  \begin{equation}\label{defphi12}
\phi_1 \ = \ \dfrac{\nu\omega_2+\frac{2}{3}-\frac{\beta b}{3}}{\epsilon(\nu-\frac{2\beta b}{3}+\frac{\beta^2b^2}{3})}-\dfrac{\omega_1}{\epsilon} ,
  \end{equation}
with $\omega_1$ and $\omega_2$ solutions to the below first order linear differential equations respectively
   \begin{equation}\label{defomega12}
\omega_1' \ + \  \left(\dfrac{1}{\beta b}+\dfrac{1-\beta b}{(1-\beta b)^2-\frac{3}{\bo}}\right)\omega_1 = \dfrac{\beta b-1}{\beta b\big((1-\beta b)^2-\frac{3}{\bo}\big)}, \end{equation}
\begin{equation}\label{defomega22}
\omega_2' \ + \  \left(\dfrac{1}{\beta b}+\dfrac{3(1-\beta b)}{(1-\beta b)^2-\frac{3}{\bo}}\right)\omega_2 = \dfrac{3(\beta b-1)}{\beta b\big((1-\beta b)^2-\frac{3}{\bo}\big)}, \end{equation} 
and $\varsigma=1$. Some additional key restrictions are required for the validity of our model~\eqref{eq:Serre1}, namely $\nu \neq 0$, $\beta b \neq 0$ and $\beta b \neq 1 \pm \sqrt{\dfrac{3}{\bo}}$. Those restrictions correspond to~\eqref{condbot} with $\gamma=0$ and $\delta=1$. Moreover, it is clear now that the r.h.s terms of the second equation of~\eqref{eq:Serre1} denoted by $\mathcal{A}$, $\mathcal{B}$, $\mathcal{C}$, $\mathcal{D}$, $\mathcal{E}$, and $\mathcal{F}$ (see Appendix~\ref{appendixB}) become (after setting $\gamma=0$ and $\delta=1$):
           \begin{equation}\label{defA}
      \mathcal{A}=-\epsilon\beta(1+\beta\omega_1 b)(1-\beta b)\partial_x^2b,
           \end{equation}
            \begin{equation}\label{defB}
      \mathcal{B}=\epsilon\beta(1+\beta\omega_1 b)(1-\beta b)\partial_x b,
      \end{equation}
        \begin{equation}\label{defC}
      \mathcal{C}=0,
      \end{equation}
      \begin{equation}\label{defD}
      \mathcal{D}=-\dfrac{2}{3}\epsilon(1+\beta\omega_1 b)(1-\beta b)^2,
\end{equation}
\begin{equation}\label{defE}
      \mathcal{E}=-\dfrac{\epsilon\beta}{2}(1+\beta\omega_1 b)(1-\beta b)\partial_x^3b,
      \end{equation}
      \begin{equation}\label{defF}
      \mathcal{F}=\frac{\beta}{2} q_1(\epsilon\zeta,\beta b)(1-\beta b)\partial_x^2 b+\frac{\epsilon\beta}{2}(1+\beta w_1 b)\zeta\partial_x^2 b.
       \end{equation}  
At this stage, it is worth mentioning that after multiplying its second equation by $\dfrac{h}{q_1(\epsilon\zeta,\beta b)}$, system~\eqref{eq:Serre1} is equivalent to the system studied in~\cite{Haidar2018} when all terms of order $\OO(\mu^2,\mu\epsilon^2)$ are ignored in the latter. This equivalence is illustrated in Appendix~\ref{appendixA}. Consequently, in the absence of surface tension (i.e. setting $\bo^{-1}=0$), system~\eqref{eq:Serre1} is also equivalent to the system studied in~\cite{Israwi11} when all terms of order $\OO(\mu^2,\mu\epsilon^2)$ are ignored in the latter.
 
The results of the previous section apply as a particular case,  allowing us to conclude that the GNCH with uneven bottoms in the water-wave case taking into account surface tension effects~\eqref{eq:Serre1} is fully justified ( i.e consistent, well-posed and convergent).
Treating the capillary term $\dfrac{\mu}{\bo}\partial_x^3\zeta$ (thanks to the suitable choice of $\nu$ in~\eqref{defnu2}) allowed us to provide a well-posedness result in the standard Hyperbolic space $H^{s}(\RR)\times H^{s+1}(\RR)$, with $s>3/2$. Unlike in \cite{Haidar2018},  controlling the capillary term necessitates the definition of the energy norm $\vert(\zeta,v)\vert^2_{X^s}\equiv \vert \zeta\vert^2 _{H^s}+\dfrac{\mu}{\bo}\vert \partial_x \zeta\vert^2 _{H^s}+\vert v\vert^2 _{H^s}+ \mu\vert \partial_xv\vert^2 _{H^s}$, and thus requests more regularity on $\zeta$. Indeed, the second term in $\vert \cdot\vert_{X^s}$ is missing from the natural energy adopted in our work $\vert \cdot \vert_{Y^s}$ (see Definition~\ref{defispace}). This term is critical in controlling the capillary term that appears in the original model studied  \cite{Haidar2018},  which yields  the well-posedness in smoother Hyperbolic space $H^{s+1}(\RR)\times H^{s+1}(\RR)$, with $s>3/2$.

\appendix  
\section{Equivalence with a model in the literature}\label{appendixA}
In this Appendix we show that the GNCH ($\epsilon=\OO(\sqrt{\mu})$) with uneven bottoms in the water-wave case taking into account surface tension effects~\eqref{eq:Serre1} is equivalent to the system studied in~\cite{Haidar2018} when neglecting in the latter all terms of order $\OO(\mu^2,\mu\epsilon^2)$. One can easily notice that the evolution equation on the surface deformation $\zeta$ of system~\eqref{eq:Serre1} is the same as the first equation of system (2) in~\cite{Haidar2018}. To show the equivalence between the evolution equation on the layer-mean horizontal velocity $v$ of system~\eqref{eq:Serre1} and the second equation of system (2) in~\cite{Haidar2018} we multiply the former equation by $\dfrac{h}{q_1(\epsilon \zeta,\beta b)}$. In what follows we omit to write the dependence on $\epsilon \zeta$ and $\beta b$ of $q_1$ and $q_2$.   Before starting,  we will state two approximations that will be extensively used in what follows:
\begin{equation}\label{approx1}
\dfrac{1}{q_1}= \dfrac{1}{1+\beta \omega_1 b} + \OO(\epsilon),
\end{equation}
and 
\begin{equation}\label{approx2}
\dfrac{1}{q_1}= \dfrac{1}{1+\beta (\omega_1+\epsilon \phi_1) b} \Bigg(1-\dfrac{\epsilon\kappa_1 \zeta}{1+\beta (\omega_1+\epsilon \phi_1) b}\Bigg) + \OO(\epsilon^2).
\end{equation}
In fact, lets consider first the multiplication of  ${\mathfrak T}[\epsilon\zeta,\beta b]V$ by $\dfrac{h}{q_1}$, one gets:
 \begin{equation*}
      \dfrac{h}{q_1}{\mathfrak T}[\epsilon\zeta,\beta b]V \ = \ hV \ + \ \mu\epsilon\beta h \dfrac{\kappa_0}{q_1} \partial_x \zeta \partial_x b V \ - \ \mu\nu h \dfrac{\partial_x\Big(q_2\partial_xV \Big)}{q_1}.
      \end{equation*}
Using~\eqref{approx1} and the definition of $\kappa_0$ in~\eqref{defkappa02} one gets the following approximation,
 \begin{equation}\label{hI/q1}
      \dfrac{h}{q_1}{\mathfrak T}[\epsilon\zeta,\beta b]V \ = \ hV \ + \ \mu\epsilon\beta h \partial_x \zeta \partial_x b V  \ - \ \mu h \nu \dfrac{\partial_x q_2} {q_1} \partial_x V\ - \ \mu h \nu \dfrac{q_2 }{q_1}\partial_x^2 V +\OO(\mu\epsilon^2).
   \end{equation}
Now lets consider respectively the third and fourth terms of the r.h.s of~\eqref{hI/q1}. In fact using the definition of $q_2$, one has
  \begin{equation}\label{nudxq2/q1}
      \nu \dfrac{\partial_x q_2}{q_1}= \dfrac{\nu \kappa_2}{q_1} \epsilon \partial_x  \zeta  +  \dfrac{\nu \kappa_2'}{q_1} \epsilon \beta \zeta \partial_x b + \dfrac{\nu \omega_2 + \nu\beta \omega_2 ' b}{q_1}\beta \partial_x b.\end{equation}
Using again~\eqref{approx1} and the definition of $\kappa_2$ in~\eqref{defkappa22}, the following approximation hold:
\begin{equation}\label{nukappa2/q1}
\dfrac{\nu\kappa_2}{q_1}=\dfrac{\nu\kappa_2}{1+\beta \omega_1 b}+\OO(\epsilon)=(1-\beta b) +\OO(\epsilon).
\end{equation}
Moreover, using~\eqref{defomega1} (with $\gamma=0$ and $\delta=1$) together with the definition of $\kappa_1$ in~\eqref{defkappa12} one has
\begin{equation}\label{nukappa'2/q1}
\dfrac{\nu\kappa_2'}{q_1}=\dfrac{\nu\kappa_2'}{1+\beta \omega_1 b}+\OO(\epsilon)=-2 -\dfrac{3}{\bo(1-\beta b)^2-3} +\OO(\epsilon).
\end{equation}
Now, using~\eqref{defomega2} (with $\gamma=0$ and $\delta=1$) one can check that the following hold:
\begin{align}\label{nuomega2+}
\dfrac{\nu \omega_2 + \nu \beta \omega_2' b}{q_1} &= -(1-\beta b) \dfrac{(1+\beta (\omega_1 +\epsilon \phi_1) b)}{q_1}\nn\\&=-(1-\beta b) \dfrac{(q_1 -\epsilon \kappa_1 \zeta)}{q_1}\ \ \text{(using definition of $q_1$)}\nn \\&=-(1-\beta b)+ \dfrac{(1-\beta b) \epsilon \kappa_1 \zeta}{q_1}\nn \\&=-(1-\beta b)+ \dfrac{(1-\beta b) \epsilon \kappa_1 \zeta}{1+\beta \omega_1 b} + \OO(\epsilon^2) \ \ \text{(using~\eqref{approx1})} \nn \\&=-(1-\beta b)+ \dfrac{\epsilon \zeta (1-\beta b)^2}{(1-\beta b)^2 -3/\bo} + \OO(\epsilon^2) \ \ \text{(using~\eqref{defkappa12})} \nn \\&=-(1-\beta b)+ \epsilon \zeta + \dfrac{3\epsilon \zeta}{\bo(1-\beta b)^2 -3} + \OO(\epsilon^2).
\end{align}
Thus, one can deduce using~\eqref{nukappa2/q1},~\eqref{nukappa'2/q1} and~\eqref{nuomega2+} that~\eqref{nudxq2/q1} can be approximated as follows:
  \begin{equation}\label{nudxq2/q1_approx}
      \nu \dfrac{\partial_x q_2}{q_1}= (1-\beta b) \epsilon \partial_x \zeta-h \beta \partial_x b +\OO( \epsilon^2).
     \end{equation}
Now, lets consider the fourth term of the r.h.s of~\eqref{hI/q1}. Multiplying~\eqref{approx2} by $q_2$, the following approximation easily follow:
\begin{equation}\label{q2/q1}
\dfrac{q_2}{q_1}= \dfrac{1+\beta \omega_2 b}{1+\beta (\omega_1+\epsilon \phi_1) b}+\dfrac{\epsilon \kappa_2 \zeta}{1+\beta \omega_1 b}-\dfrac{\epsilon\kappa_1 \zeta(1+\beta \omega_2 b)}{(1+\beta (\omega_1+\epsilon \phi_1) b)^2} + \OO(\epsilon^2).
\end{equation}
Using~\eqref{defphi12}, the first term of the r.h.s of~\eqref{q2/q1} can be written as:
\begin{equation}\label{1-q2/q1}
 \dfrac{1+\beta \omega_2 b}{1+\beta (\omega_1+\epsilon \phi_1) b}=\dfrac{\bo(1-\beta b)^2-3}{\bo-3}.
\end{equation}
Using~\eqref{defkappa22}, the second term of the r.h.s of~\eqref{q2/q1} can be written as:
\begin{equation}\label{2-q2/q1}
 \dfrac{\epsilon \kappa_2 \zeta}{1+\beta \omega_1 b}=\dfrac{\epsilon \zeta (1-\beta b)}
 {\nu}.
\end{equation}
Using the definition of $\kappa_1$ in~\eqref{defkappa12} and~\eqref{1-q2/q1}, the third term of the r.h.s of~\eqref{q2/q1} can be approximated as follows:
\begin{align}\label{3-q2/q1}
\dfrac{\epsilon\kappa_1 \zeta(1+\beta \omega_2 b)}{(1+\beta (\omega_1+\epsilon \phi_1) b)^2}&= \dfrac{\epsilon \zeta(1+\beta \omega_1 b) (1-\beta b)\bo}{(1+\beta (\omega_1+\epsilon \phi_1)b)(\bo-3)}\nn\\&= \dfrac{\epsilon \zeta (1-\beta b)\bo}{(\bo-3)} + \OO(\epsilon^2).
\end{align}
Thus, one can deduce using~\eqref{1-q2/q1},~\eqref{2-q2/q1} and~\eqref{3-q2/q1} that~\eqref{q2/q1} can be approximated as follows:
\begin{equation}\label{q2/q1_approx}
\dfrac{q_2}{q_1}= \dfrac{\bo(1-\beta b)^2-3}{\bo-3} +\dfrac{\epsilon \zeta (1-\beta b)}
 {\nu}-\dfrac{\epsilon \zeta (1-\beta b)\bo}{(\bo-3)} + \OO(\epsilon^2).
\end{equation}
Now, using~\eqref{nudxq2/q1_approx} and~\eqref{q2/q1_approx}, one can deduce that approximation~\eqref{hI/q1} can be written as follows:
 \begin{align}\label{hI/q1_approx}
      \dfrac{h}{q_1}{\mathfrak T}[\epsilon\zeta,\beta b]V &=  hV \ + \ \mu\epsilon\beta h \partial_x \zeta \partial_x b V  \ - \ \mu h\Big((1-\beta b) \epsilon \partial_x \zeta -h \beta \partial_x b \Big) \partial_x V\nn \\ & \qquad -  \mu h\Big(\dfrac{(1-\beta b)^2}{3}-\dfrac{1}{\bo}+\dfrac{2}{3}\epsilon \zeta (1-\beta b)      \Big)\partial_x^2 V +\OO(\mu\epsilon^2).\end{align}
Now lets recall the notation used at the beginning of Section 3,  $$\T[h,b]V\equiv \dfrac{-1}{3h}\partial_x(h^3\partial_x V)+ \dfrac{1}{2h}[\partial_x(h^2(\partial_x b) V)-h^2(\partial_x b)(\partial_x V)]+(\partial_x b)^2 V .$$
In fact, one can notice that $\T[h,\beta b]V$ can be approximated  as follows:
\begin{align}\label{T_approx}\T[h,\beta b]V &= -(1-\beta b) \epsilon \partial_x \zeta \partial_x V + h \beta \partial_x b \partial_x V -\dfrac{1}{3}(1-\beta b)^2 \partial_x^2 V -\dfrac{2}{3} \epsilon \zeta(1-\beta b)\partial_x^2 V \nn \\ & \quad + \beta \epsilon \partial_x \zeta \partial_x b V +\dfrac{\beta}{2} h \partial_x^2 b V + \OO(\epsilon^2).\end{align}
Thus, using~\eqref{T_approx}, the approximation~\eqref{hI/q1_approx} can be rewritten as:
    \begin{equation}\label{hI/q1_approx2}
      \dfrac{h}{q_1}{\mathfrak T}[\epsilon\zeta,\beta b]V =(h+\mu h \T [h,\beta b])V-\dfrac{\mu\beta}{2}h^2 \partial_x^2 b V + \dfrac{\mu h}{\bo} \partial_x^2 V+\OO(\mu\epsilon^2). 
   \end{equation}
Now, using~\eqref{approx1} and the definition of $\mathcal{F}$ in~\eqref{defF}, one can check that the following approximation holds:

     \begin{equation}\label{F_approx}
  \mu \dfrac{h}{q_1} [\mathcal{F}]=  \dfrac{\mu\beta}{2} h^2 \partial_x^2 b +\OO(\mu\epsilon^2), 
   \end{equation} 
thus,~\eqref{hI/q1_approx2} can be rewritten as:
     \begin{equation}\label{hI/q1_approx3}
      \dfrac{h}{q_1}{\mathfrak T}[\epsilon\zeta,\beta b]V =(h+\mu h \T [h,\beta b])V- \mu \dfrac{h}{q_1} [\mathcal{F}] V + \dfrac{\mu h}{\bo} \partial_x^2 V+\OO(\mu\epsilon^2). 
   \end{equation}
At this stage, one can notice that after multiplying  the second equation of~\eqref{eq:Serre1} by $\dfrac{h}{q_1}$, and using the \emph{BBM} trick represented by $\partial_t v + \epsilon v \partial_x v = -\partial_x \zeta +\OO(\mu)$, one gets the following approximation:
\begin{align}\label{eq2_lhs:Serre1} \dfrac{h}{q_1}\mathfrak T[\epsilon\zeta,\beta b] \left( \partial_{ t}   v + \epsilon { v } \partial_x {  v } \right) +h \partial_x \zeta -\mu\dfrac{h}{q_1} [\mathcal{F}] \partial_x \zeta \nn \\ =(h+\mu h \T [h,\beta b])(\partial_{ t}   v + \epsilon { v } \partial_x {  v })+ h \partial_x \zeta - \dfrac{\mu h}{\bo} \partial_x^3 \zeta+\OO(\mu\epsilon^2). \end{align}
Moreover, using the definitions of $\mathcal{A}$, $\mathcal{B}$, $\mathcal{C}$, $\mathcal{D}$ and $\mathcal{E}$ in~\eqref{defA},~\eqref{defB},~\eqref{defC},~\eqref{defD} and~\eqref{defE} respectively and using~\eqref{approx1} one has:
\begin{align} \label{eq2_rhs:Serre1}
& \dfrac{h}{q_1}\cdot\Big[ \mu\Big([\mathcal{A}]v\partial_x v +  [\mathcal{B}](\partial_x v)^2 +[\mathcal{D}]\partial_x\big((\partial_x v)^2\big)+[\mathcal{E}]v^2\Big) \Big]    \\ 
&=   -\mu\epsilon h\big[\beta (1-\beta b)\partial_x^2 b v\partial_x v -\beta\partial_x b (1-\beta b) (\partial_x v)^2 + \frac{2}{3} (1-\beta b)^2 \partial_x ((\partial_x v)^2)+\frac{\beta}{2} (1-\beta b)\partial_x^3 b v^2\big]+\OO(\mu\epsilon^2)\nonumber  \\
& =-\mu\epsilon h \Q[h,\beta b]v +\OO(\mu\epsilon^2) \nonumber  \;,
 \end{align}     
with $ \Q[h,\beta b]v=\dfrac{2}{3h}\partial_x(h^3 (\partial_x v)^2)+\beta h (\partial_x v)^2 \partial_x b +\dfrac{\beta}{2h} \partial_x(h^2 v^2 \partial_x^2 b) +\beta^2 \partial_x^2 b\partial_x b v^2.$\\
\\
Finally, using~\eqref{eq2_lhs:Serre1} and~\eqref{eq2_rhs:Serre1}, it becomes clear now that the second equation of system~\eqref{eq:Serre1} multiplied by $\dfrac{h}{q_1}$ is equivalent to the following equation:
\begin{equation*}
(h+\mu h \T [h,\beta b])(\partial_{ t}   v + \epsilon { v } \partial_x {  v })+ h (1-\dfrac{\mu}{\bo}\partial_x^2)\partial_x \zeta =-\mu\epsilon h \Q[h,\beta b]v+\OO(\mu\epsilon^2).
\end{equation*}
The above equation corresponds exactly to the second equation of system (2) in~\cite{Haidar2018} when neglecting in the latter all terms of order $\OO(\mu^2,\mu\epsilon^2)$. Consequently, in the absence of surface tension (i.e. setting $\bo^{-1}=0$) and after multiplying its second equation by $\dfrac{h}{q_1}$, system~\eqref{eq:Serre1} is also equivalent to the system studied in~\cite{Israwi11} when neglecting in the latter all terms of order $\OO(\mu^2,\mu\epsilon^2)$.

\section{Functions}\label{appendixB}
We conclude with a detailed statement of the functions $\mathcal{A}$, $\mathcal{B}$, $\mathcal{C}$, $\mathcal{D}$, $\mathcal{E}$, and $\mathcal{F}$ given in Section~\ref{section3}.

 \begin{eqnarray*} \mathcal{A}&=&2\epsilon\beta(1+\beta w_2 b) \partial_x b \nu'(\beta b)  \partial_x(f(\beta b)^2-\gamma g(\beta b)^2)+2\epsilon\beta^2 \nu \partial_x b w'_2 b  \partial_x(f(\beta b)^2-\gamma g(\beta b)^2)\nn\\&\qquad+& 2\epsilon\beta w_2 \nu \partial_x b  \partial_x(f(\beta b)^2-\gamma g(\beta b)^2)+3\epsilon\beta w_2 \nu b  \partial_{x}^2(f(\beta b)^2-\gamma g(\beta b)^2)-\dfrac{3\epsilon}{\bo}\partial_{x}^2(f(\beta b)^2-\gamma g(\beta b)^2)\nn\\&\qquad-&3\epsilon\beta w_1 b [\lambda(\beta b)]\partial_x^2 (f(\beta b)^2-\gamma g(\beta b)^2) \nn\\&\qquad+&\epsilon(1+\beta w_1 b)\theta[2\beta\partial_x(\partial_x b g'(\beta b))+\partial_x^2 g(\beta b)]+\epsilon(1+\beta w_1 b)[2\theta+(\gamma-1)g(\beta b)]\partial_x(\beta \partial_x b g'(\beta b)) \nn\\&\qquad+&\epsilon\beta(1+\beta w_1 b) \alpha(f(\beta b)^2-\gamma g(\beta b)^2)\partial_x^2 b +4\epsilon\beta(1+\beta w_1 b)\alpha\partial_{x}(f(\beta b)^2-\gamma g(\beta b)^2)\partial_x b\nn\\&\qquad+&3\epsilon\beta(1+\beta w_1 b)[\frac\gamma3 f(\beta b)+ \frac23\delta^{-1}f(\beta b)] \partial_{x}^2(f(\beta b)^2-\gamma g(\beta b)^2)b+\epsilon\beta(1+\beta w_1 b)(\theta_1-\alpha_1)g(\beta b)\partial_x^2 b \nn\\&\qquad+&\epsilon\beta(1+\beta w_1 b)(2\theta_1-\alpha_1)\beta(\partial_x b)^2 g'(\beta b)\nn\\&\qquad+& \epsilon\beta(1+\beta w_1 b)[2\theta_1-2\alpha_1 +\frac13(\delta^{-1} -\beta b)f'(\beta b)-\frac\gamma3g'(\beta b)](\beta\partial_x b g'(\beta b) +\partial_x(g(\beta b)))\partial_x b\nn\\&\qquad+& \epsilon \beta^2(1+\beta w_1 b)\eta(\beta b)(\partial_x b)^2(f(\beta b)^2-\gamma g(\beta b)^2)+\epsilon\beta^2(1+\beta w_1 b) \frac\gamma3 f'(\beta b) b\partial_x ^2 b(f(\beta b)^2-\gamma g(\beta b)^2)\nn\\&\qquad+& 2\epsilon\beta^2(1+\beta w_1 b)[\frac{2\gamma}{3} f'(\beta b)] b\partial_x b \partial_x (f(\beta b)^2-\gamma g(\beta b)^2) -\epsilon\beta^2(1+\beta w_1 b) f(\beta b)\partial_x ^2 (f(\beta b)^2-\gamma g(\beta b)^2)\nn\\&\qquad+&\epsilon \beta^2(1+\beta w_1 b) \eta_1(\beta b)(\partial_x b)^2 g(\beta b)+\epsilon\beta^3(1+\beta w_1 b)[\frac\gamma 3f''(\beta b)](\partial_x b)^2 b(f(\beta b)^2-\gamma g(\beta b)^2)\nn\\ &\qquad -& \beta \partial_x b \nu'(\beta b)(1+\beta w_2 b) \epsilon \partial_x (\varsigma)-\epsilon \nu(1+\beta w_2 b)\partial_x^2 (\varsigma)\nn\\&\qquad-&  \nu\beta [\partial_x (w_2) b +w_2 \partial_x b]\epsilon \partial_x (\varsigma)+\epsilon(1+\beta w_1 b)[2 s(\beta b) +\partial_x (t(\beta b))].\end{eqnarray*}
\begin{eqnarray*}\mathcal{B}&=&\epsilon\beta \partial_x b \nu'(\beta b)(1+\beta w_2 b)(f(\beta b)^2-\gamma g(\beta b)^2)+\epsilon\beta^2 \nu \partial_x b w_2'(\beta b) b (f(\beta b)^2-\gamma g(\beta b)^2)\nn\\&\qquad+&\epsilon\beta \nu w_2\partial_x b (f(\beta b)^2-\gamma g(\beta b)^2)+ 3\epsilon\beta w_2 b \nu \partial_x (f(\beta b)^2-\gamma g(\beta b)^2)-\dfrac{3\epsilon}{\bo}\partial_x (f(\beta b)^2-\gamma g(\beta b)^2))\nn\\&\qquad-&3\epsilon\beta w_1 b [\lambda(\beta b)]\partial_x (f(\beta b)^2-\gamma g(\beta b)^2)\nn\\&\qquad+&\epsilon(1+\beta w_1 b)[2\theta+(\gamma-1)g(\beta b)](\beta\partial_x b g'(\beta b)+\partial_x (g(\beta b)))\nn\\&\qquad+&\epsilon\beta(1+\beta w_1 b)[2\alpha](f(\beta b)^2-\gamma g(\beta b)^2)\partial_x b\nn\\&\qquad+& 3\epsilon\beta(1+\beta w_1 b)[\frac\gamma3 f(\beta b)+ \frac23\delta^{-1}f(\beta b)]\partial_x (f(\beta b)^2-\gamma g(\beta b)^2) b\nn\\&\qquad+&\epsilon\beta(1+\beta w_1 b)(2\theta_1 -\alpha_1)g(\beta b)\partial_x b+\epsilon\beta^2(1+\beta w_1 b)[\frac{2\gamma}{3}f'(\beta b)] b\partial_x b (f(\beta b)^2-\gamma g(\beta b)^2) \nn\\&\qquad-&\epsilon\beta^2(1+\beta w_1 b)  f(\beta b)  b^2 \partial_x (f(\beta b)^2-\gamma g(\beta b)^2)-\epsilon\beta (1+\beta w_2 b)\partial_x b \nu'(\beta b) \varsigma-2\epsilon\nu(1+\beta w_2 b) \partial_x(\varsigma)\nn\\&\qquad-&\epsilon\beta \nu[\partial_x (w_2) b +w_2 \partial_x b]\varsigma+\epsilon(1+\beta w_1 b)[\frac12 \partial_x((1-\gamma)g(\beta b)^2)+t(\beta b)].\end{eqnarray*}
\begin{eqnarray*}\mathcal{C}&=&\epsilon\beta \partial_x b \nu'(\beta b)(1+\beta w_2 b)(f(\beta b)^2-\gamma g(\beta b)^2)+\epsilon\beta^2 \nu \partial_x b w_2'(\beta b) b (f(\beta b)^2-\gamma g(\beta b)^2)\nn\\&\quad+&\epsilon\beta \nu w_2\partial_x b (f(\beta b)^2-\gamma g(\beta b)^2)+ 3\epsilon\beta w_2 b \nu \partial_x (f(\beta b)^2-\gamma g(\beta b)^2)-\dfrac{3\epsilon}{\bo}\partial_x (f(\beta b)^2-\gamma g(\beta b)^2))\nn\\&\quad-&3\epsilon\beta w_1 b [\lambda(\beta b)]\partial_x (f(\beta b)^2-\gamma g(\beta b)^2)\nn\\&\quad+&\epsilon(1+\beta w_1 b) [\theta] (\beta\partial_x b g'(\beta b) + 2\partial_x (g(\beta b)))+\epsilon(1+\beta w_1 b)  [\theta + \frac23 (\gamma-1) g(\beta b)] \beta \partial_x b g'(\beta b)\nn\\&\qquad+& \epsilon\beta(1+\beta w_1 b) [2\alpha](f(\beta b)^2-\gamma g(\beta b)^2)\partial_x b \nn\\&\qquad+&3\epsilon\beta(1+\beta w_1 b)[\frac\gamma3 f(\beta b)+ \frac23\delta^{-1}f(\beta b)]\partial_x (f(\beta b)^2-\gamma g(\beta b)^2) b\nn\\&\qquad+& \epsilon\beta (1+\beta\omega_1 b)[2\theta_1-2\alpha_1 +\frac13 (\delta^{-1}-\beta b)f'(\beta b)-\frac\gamma3 g'(\beta b)]g(\beta b)\partial_x b\nn\\&\qquad+& \epsilon\beta^2(1+\beta w_1 b)[\frac{2\gamma}{3}f'(\beta b)] b\partial_x b (f(\beta b)^2-\gamma g(\beta b)^2)-\epsilon\beta^2(1+\beta w_1 b) f(\beta b)  b^2 \partial_x (f(\beta b)^2-\gamma g(\beta b)^2)\nn\\&\qquad-&\epsilon\beta \partial_x b \nu'(\beta b)(1+\beta w_2 b)\varsigma-2\epsilon\nu(1+\beta w_2 b) \partial_x (\varsigma)  -\epsilon\beta \nu[\partial_x (w_2) b +w_2 \partial_x b]\varsigma\nn\\&\qquad+& \epsilon(1+\beta w_1 b)[\frac13 \partial_x ( (1-\gamma)g(\beta b)^2))+t(\beta b)].\end{eqnarray*}
\begin{eqnarray*}\mathcal{D}&=&3\frac{\epsilon}{2}\nu\beta w_2 b (f(\beta b)^2-\gamma g(\beta b)^2) -\frac{3\epsilon}{2\bo}(f(\beta b)^2-\gamma g(\beta b)^2))-\frac{3\epsilon}{2}\beta w_1 b [\lambda(\beta b)](f(\beta b)^2-\gamma g(\beta b)^2)\nn\\&\qquad+&\frac{\epsilon}{2}(1+\beta w_1 b)[2\theta+(\gamma-1)g(\beta b)]g(\beta b)\nn\\&\qquad+&\frac{\epsilon}{2}(1+\beta w_1 b)[\theta+\frac23(\gamma-1)g(\beta b)]g(\beta b)+3\frac{\epsilon\beta}{2}(1+\beta w_1 b)[\frac\gamma3 f(\beta b)+ \frac23\delta^{-1}f(\beta b)](f(\beta b)^2-\gamma g(\beta b)^2)b\nn\\&\qquad-& 3\frac{\epsilon\beta^2}{2}(1+\beta w_1 b)[\frac13 f(\beta b)]b^2(f(\beta b)^2-\gamma g(\beta b)^2))-\frac{3  \nu \epsilon \varsigma}{2}(1+\beta w_2 b)+\frac{2   \epsilon (1-\gamma)g(\beta b)^2}{3}(1+\beta w_1 b)\nn\\&=&\dfrac{2\epsilon}{3}(1+\beta\omega_1 b)(\gamma-1)g(\beta b)^2.
\end{eqnarray*}
\begin{eqnarray*}\mathcal{E}&=&\frac{\epsilon}{2}\beta \partial_x b \nu'(\beta b)(1+\beta w_2 b) \partial_x^2 (f(\beta b)^2-\gamma g(\beta b)^2)+\frac{\epsilon}{2}\beta^2\nu\partial_x b w_2' b \partial_x^2 (f(\beta b)^2-\gamma g(\beta b)^2)\nn\\&\qquad+&\frac{\epsilon}{2}\nu \beta w_2 \partial_x b \partial_x^2 (f(\beta b)^2-\gamma g(\beta b)^2))+\frac{\epsilon}{2}\nu \beta w_2 b  \partial_x^3 (f(\beta b)^2-\gamma g(\beta b)^2)-\frac{\epsilon}{2\bo}\partial_x^3 (f(\beta b)^2-\gamma g(\beta b)^2)\nn\\&\qquad-&\frac{\epsilon}{2}\beta w_1 b \lambda(\beta b) \partial_x^3(f(\beta b)^2-\gamma g(\beta b)^2) \nn\\&\qquad+&\epsilon \beta (1+\beta w_1 b)[\theta]\partial_x^2(\partial_x b g'(\beta b))+\frac{\epsilon}{2} \beta (1+\beta w_1 b)[\alpha]\partial_x(f(\beta b)^2-\gamma g(\beta b)^2)\partial_x^2 b\nn\\&\qquad+&\frac{\epsilon}{2} \beta (1+\beta w_1 b)[2\alpha]\partial_x^2(f(\beta b)^2-\gamma g(\beta b)^2)\partial_x b\nn\\&\qquad+&\frac{\epsilon}{2} \beta (1+\beta w_1 b)[\frac\gamma3 f(\beta b)+ \frac23\delta^{-1}f(\beta b)]\partial_x^3(f(\beta b)^2-\gamma g(\beta b)^2) b\nn\\&\qquad+&
      \epsilon\beta^2(1+\beta w_1 b)[\theta_1 -\alpha_1]\partial_x b g'(\beta b)\partial_x^2 b
      \nn\\&\qquad+&\epsilon\beta(1+\beta w_1 b)[2\theta_1-2\alpha_1 +\frac13(\delta^{-1} -\beta b)f'(\beta b)-\frac\gamma3g'(\beta b)]\partial_x(\beta\partial_x b g'(\beta b))\partial_x b \nn\\&\qquad+&
      \frac{\epsilon}{2}\beta^2 (1+\beta w_1 b)[\eta](\partial_x b)^2\partial_x (f(\beta b)^2-\gamma g(\beta b)^2)
      +\frac{\epsilon}{2}\beta^2 (1+\beta w_1 b)[\frac\gamma3 f'(\beta b)]b\partial_x^2 b\partial_x (f(\beta b)^2-\gamma g(\beta b)^2)\nn\\&\qquad+&\frac{\epsilon}{2}\beta^2 (1+\beta w_1 b)[\frac{2\gamma}{3} f'(\beta b)]b\partial_x b\partial_x^2 (f(\beta b)^2-\gamma g(\beta b)^2)
      -\frac{\epsilon}{2}\beta^2 (1+\beta w_1 b)[\frac13 f(\beta b)]b^2\partial_x^3 (f(\beta b)^2-\gamma g(\beta b)^2)\nn\\&\qquad+&\epsilon\beta^2 (1+\beta w_1 b)[\eta_1](\partial_x b)^2 \beta\partial_x b g'(\beta b)+\frac{\epsilon}{2}\beta^3(1+\beta w_1 b)[\frac{\gamma}{3} f''(\beta b)](\partial_x b)^2 b \partial_x(f(\beta b)^2-\gamma g(\beta b)^2)\nn\\&\qquad+&\epsilon(1+\beta w_1 b) \partial_x(s(\beta b)).\end{eqnarray*}
\begin{eqnarray*}
\mathcal{F}&=&\beta q_1(\epsilon\zeta,\beta b)[\alpha](\gamma+\delta)\partial_x^2 b+\epsilon\beta(1+\beta w_1 b)(\gamma+\delta)[\theta_1-\alpha_1]\partial_x^2 b \zeta+\beta^2q_1(\epsilon\zeta,\beta b)[\eta](\gamma+\delta)(\partial_x b)^2\nn\\&\qquad+&\beta^2q_1(\epsilon\zeta,\beta b)[\frac\gamma3 f'(\beta b)](\gamma+\delta)b\partial_x^2 b+\epsilon\beta^2(1+\beta w_1 b)(\gamma+\delta)[\eta_1](\partial_x b)^2 \zeta\nn\\&\qquad+&\beta^3q_1(\epsilon\zeta,\beta b)[\frac\gamma3 f''(\beta b)](\gamma+\delta)b(\partial_x b)^2.\end{eqnarray*}

\end{document}